\RequirePackage{fix-cm}
\documentclass{svjour3}                     
\smartqed  
\usepackage{graphicx}
 \usepackage{mathptmx}      
\usepackage[english]{babel}
\usepackage{epsfig}
\usepackage{amsmath}
\usepackage{amssymb}
\usepackage{amsbsy}
\usepackage{amsfonts}
\usepackage{mathrsfs} 
\usepackage{array}
\usepackage{verbatim}
\usepackage{graphicx}
\usepackage[dvipsnames]{xcolor}
\usepackage[noadjust]{cite}
\usepackage{tabularx}
\usepackage{multirow}
\usepackage{bm}
\usepackage{booktabs}
\usepackage{subfig}
\usepackage{stmaryrd}
\usepackage{textcomp}
\usepackage[font=footnotesize]{caption}
\newcommand{\revision}[1]{{\color{black}{#1}}}

\journalname{Journal of Scientific Computing}
\begin{document}

\title{A seamless, extended DG approach for \revision{advection-diffusion} problems on unbounded domains\thanks{TB and LB were supported were supported by the ESCAPE-2 project of the Horizon 2020 research and innovation programme (grant agreement No 800897)}
}


\author{Federico Vismara*\and Tommaso Benacchio \and Luca Bonaventura}

\authorrunning{F. Vismara et al.} 

\institute{MOX -- Modelling and Scientific Computing\\ 
           Dipartimento di Matematica, Politecnico di Milano \\
           Via Bonardi 9, 20133 Milano, Italy \\\vspace{1mm}
           \\
           * Corresponding author. \\Present address: CASA, Technische Universiteit Eindhoven, The Netherlands.\\\email{f.vismara@tue.nl}}

\date{Received: date / Accepted: date}

\maketitle

\begin{abstract}
We propose and analyze a seamless extended Discontinuous Galerkin (DG) discretization of \revision{advection-diffusion} equations on semi-infinite domains. The semi-infinite half line is split into a finite subdomain where the model uses a standard polynomial basis, and a semi-unbounded subdomain where scaled Laguerre functions are employed as basis and test functions. Numerical fluxes enable the coupling at the interface between the two subdomains in the same way as standard single domain DG interelement fluxes. A novel linear analysis on the extended DG model yields \revision{unconditional stability with respect to the P\'eclet number}.
Errors due to the use of different sets of basis functions on different portions of the domain are negligible, as highlighted in numerical experiments  with the linear advection-diffusion and viscous Burgers' equations. With an added damping term on the semi-infinite subdomain, the extended framework is able to efficiently simulate absorbing boundary conditions without additional conditions at the interface. A few modes in the semi-infinite subdomain are found to suffice to deal with outgoing single wave and wave train signals \revision{more accurately than standard approaches at a given computational cost}, thus providing an appealing model for fluid flow simulations in unbounded regions.
\end{abstract}

\section{Introduction}
\label{sec:intro} \indent  
The correct modelling of evolution problems over arbitrarily large regions has a wide range of applications in computational physics and poses several still unsolved challenges. An especially relevant application area is atmospheric modelling, where the region of interest to forecasts - typically, the troposphere and lower stratosphere - should not feature spurious reflections of upwardly propagating waves generated by the computational model lid, see, e.g., \cite{bonaventura:2000, klemp:1978, klemp:1983}. At the same time, as computational resources enable raising the lid, an accurate description of upper atmosphere phenomena is of paramount practical interest and the goal of fully integrated space weather models is increasingly being discussed, see, e.g., \cite{akmaev:2011, jackson:2019}.

Approximations over arbitrarily large regions usually rely on the creation of an artificial boundary separating the region of interest from an external region. Analytical approaches, e.g., \cite{engquist:1977, dea:2011, israeli:1981}, attempt to impose conditions at the artificial boundary in order to let outgoing perturbations propagate without spurious reflections. However, these conditions can be difficult to determine and computationally expensive, and may require \emph{ad hoc} information on the outgoing signal. 

An alternative approach is provided by numerical techniques based on absorbing (or sponge) layers. These are buffer regions placed next to the artificial boundary where perturbations leaving the computational domain are damped to a prescribed external solution by an artificial dissipation term. The choice of the parameters to be employed in these regions, however, is non-trivial and the corresponding buffers may be quite large and entail substantial computational costs.
A complete review of the proposed approaches for open boundary conditions is beyond the scope of this paper. Comprehensive reviews can be found, for example, in  \cite{appelo:2009, astley:2000, benacchio:2013, black:1998, gerdes:2000, rasch:1986}.

In \cite{benacchio:2013, benacchio:2019}, a numerical approach to open boundary conditions was proposed, based on the use of scaled Laguerre functions \cite{shen:2001, shen:2009b, shen:2009, wang:2009, zhuang:2010} for spectral approximations on the semi-infinite line, coupled to finite volume or finite element discretizations on a finite domain. Different approaches were used on either side of the artificial boundary, and only the hyperbolic case was considered. In addition, the stability analysis in \cite{benacchio:2019} only concerned the discretization on the semi-infinite subdomain. 

This paper \revision{develops} and completes the previously proposed approach by presenting two major improvements. First, the method is extended to conservation laws \revision{with viscous terms}, thus showing that all relevant terms of standard computational fluid dynamics models are amenable to a discretization based on scaled Laguerre functions. Second, the \revision{advection-diffusion} model equations are discretized in a completely seamless way, using a discontinuous Galerkin (DG) finite element formulation that relies on scaled Laguerre functions as both basis and test functions on the semi-infinite portion, along with Gauss-Laguerre-Radau quadrature rules for numerical integration. 

The resulting extended DG approach and its numerical stability on the entire half line are analyzed in the case of a linear advection-diffusion equation, by considering several options for the polynomial basis and quadrature rules. Scaled Laguerre functions and Gauss-Laguerre-Radau quadrature formulae prove to be the most stable in all hydrodynamic regimes. 

Numerical validation of the proposed method is carried out for the linear advection-diffusion and the viscous Burgers' equation. First, a relatively large number of basis functions are used in the semi-infinite portion of the domain, in order to evaluate the errors introduced by considering different bases on either side of the finite/semi-infinite interface. By taking as reference a standard single-domain DG discretization, spurious reflections are found to be of negligible entity. A reaction damping term is then introduced in the semi-infinite layer, in order to simulate the propagation of an isolated Gaussian initial datum and boundary generated wave train from the finite subdomain into the semi-infinite subdomain, where they are damped using an appropriate of Laguerre basis functions. Few basis functions suffice to absorb outgoing signals with minimal reflections and low computational cost, thereby providing an efficient implementation of an absorbing layer \revision{compared with standard approaches.}  

\revision{While other multi-domain strategies for the simulation of fluid dynamics in unbounded domains are already available in the literature, the analysis in the paper \cite{zhuang:2010} was focused on Stokes' equations, while the numerical examples with the Navier-Stokes equations did not include a damping term and were run for relatively short final times. Other authors  \cite{zhuang:2010b,tatari:2014} considered the discretization of even-order operators with coupled spectral-spectral discretizations in the stationary \cite{zhuang:2010b} and time-dependent \cite{tatari:2014} case -  though \cite{tatari:2014} focused on numerical validation rather than on stability analyses. In addition, the authors of \cite{gu:2021}  recently developed discretizations of elliptic problems in unbounded domains with orthogonal Jacobi rational functions. 
}

The outline of the paper is as follows. Section \ref{sec:modeq} contains the model equation and outlines the numerical discretization. Stability properties of the proposed method are analyzed in detail in Section \ref{sec:analysis}, and Section \ref{sec:experiments} contains the results of the numerical experiments.  The final Section \ref{sec:conclu} draws conclusions, discussing possible extensions and future work. The Appendix \ref{sec:alternatives} summarizes the analysis of different possible discretizations on the semi-infinite domain, based on either scaled Laguerre functions or polynomials or on different choices for the numerical quadrature rules.  As already shown in \cite{benacchio:2019}, some of these alternatives are shown to be more problematic in the advection-dominated case.

\section{The extended DG discretization approach}
\label{sec:modeq} \indent  
We consider as a model problem the one-dimensional nonlinear conservation law 
\revision{with diffusive terms} for the unknown $c$
\begin{equation}\label{eq:advdiffreacgen}
    \frac{\partial c}{\partial t} + \frac{\partial f(c)}{\partial z} = \frac{\partial}{\partial z}\left(\mu(z,t)\frac{\partial c}{\partial z}\right) + s(c,z,t)  
\end{equation}
for $z\in[0,+\infty) $ and $t\in[0, T]. $ We assume that the diffusion coefficient $\mu$ is a smooth function of its variables and that there are two positive constants $\mu_0$ and $\mu_1$ such that
\begin{equation}
    0<\mu_0\leq\mu(z,t)\leq\mu_1 \qquad \forall z\in[0,+\infty), \quad \forall t>0.
\end{equation}
For simplicity, we only consider here Dirichlet boundary conditions  
\begin{equation}\label{eq:dir_bc}
   c(0,t)=g_0(t)  \ \ \ \ \lim_{z\rightarrow +\infty} c(z,t)=0.
\end{equation}
We refer to \cite{vismara:2020} for a comprehensive analysis of other boundary conditions in the linear advection-diffusion case.
Next, drawing on the approach of \cite{benacchio:2013,benacchio:2019} we split the domain as ${\mathbb R}^+=[0,L]\cup[L,+\infty)$
and introduce an extended DG finite element discretization on ${\mathbb R}^+$ using a standard polynomial basis on $[0,L] $ and the scaled Laguerre functions as both basis and test functions on  $[L,+\infty)$. 
 More specifically, on the $[0,L]$ interval a mesh of $N$ non-overlapping elements $K_m$ of size $\Delta z_m  \leq h $ is considered, such that $[0,L]=\bigcup_{m=1}^{N} K_m$.
The center of the generic element $K_m$ is denoted by $z_m$, while $z_{m\pm1/2} $ denote its boundary points. 
The affine local maps $z = Z_{m}(\xi)= \xi\Delta z_m/2+z_{m}$ map the master element
$\hat K = [-1,1]$ onto each $K_m$.  
For each non-negative integer $p$, we then denote by $\mathbb{P}_p$ the set of all polynomials of degree less or equal to $p$ on $\hat K$. We also define
$\mathbb{P}_p(K_m)= \left\{w : w= v\circ {Z}^{-1}_{m}, \quad  v\in \mathbb{P}_p \right\}$.
For each polynomial degree $p$, the discontinuous finite element spaces are defined as:
\begin{equation}
V^p_{h}=\left\{ v \in L^2( [0,L] ) :
v|_{K_m} \in \mathbb{P}_{p}(K_m)\quad  m=1,\dots,N  \right\}.
\label{eq:fespace_disc}
\end{equation}
The bases of $\mathbb{P}_{p}(K_m)$ are obtained from Legendre polynomials as follows.
First, for $\xi \in \hat{K}$, Legendre polynomials are defined by the recurrence relation:
\begin{align}
L_{k+1} &= \frac{2k+1}{k+1} \xi L_k(\xi) - \frac{k}{k+1} L_{k-1}(\xi), & k = 1,2,\ldots \\
L_0(\xi) &= 1, \quad L_1(\xi) = \xi. 
\end{align}
Legendre polynomials form an orthogonal basis for polynomials on $\hat{K}$ since
\begin{equation}
\int_{-1}^{1} L_k(\xi) L_l(\xi) d\xi = \frac{2}{2k+1} \delta_{kl}.
\end{equation}
For each element $K_m,\,m=1,\dots,N$  we then denote by $\phi^{m}_j(z), j=0,\dots, p$  the basis and test functions given by
\begin{equation}
\phi^m_l(z)= \sqrt{2l+1}L_l\Big(2\frac{z-z_m}{\Delta z_m}\Big).
\end{equation}
Notice that the normalization is chosen so that
\begin{equation}
\label{eq:normalize_phi}
\int_{z_{m-\frac 12}}^{z_{m+\frac 12}} \phi^m_{k}(z) \phi^m_{l}(z) \,dz = \Delta z_m \delta_{kl}.
\end{equation}
Therefore, the solution of \eqref{eq:advdiffreacgen} will be represented on each subinterval $K_m$ as
\begin{equation}
\label{eq:modalexpdg}
 c(z,t)\approx \sum_{j=0}^{p}  c^{(j)}_{m}(t) \phi^m_j(z), \ \ \ \ z\in K_m
\end{equation}
and standard Gauss-Legendre formulae will be used to discretize the resulting integrals.
For the semi-infinite interval  $K_{\infty}=[L,+\infty)$, we consider the scaled Laguerre functions as modal basis. The possible alternatives are discussed and analyzed  in \cite{benacchio:2019} for the purely hyperbolic case and  in Appendix \ref{sec:alternatives} of this paper for the hyperbolic-parabolic case.  
More specifically, defining scaled Laguerre polynomials on $[0,+\infty) $ by
\begin{align}
\label{eq:laguerre_pol}
     (k+1)\mathscr{L}^\beta_{k+1}(x)&=(2k+1-\beta x)\mathscr{L}^\beta_k(x)-k\mathscr{L}^\beta_{k-1}(x), \\ 
     \mathscr{L}^\beta_0(x) &= 1,\; \mathscr{L}^\beta_1(x) = 1-\beta x, 
\end{align}
scaled Laguerre functions are defined by
\begin{equation}
    \hat{\mathscr{L}}_k^{\beta}(x)=e^{-\beta x/2}\mathscr{L}_k^{\beta}(x).
\end{equation}
for the scaling factor $\beta>0$, and are a complete orthogonal system in $L^2(\mathbb{R}^+)$, such that 
\begin{equation}
    \int_0^{+\infty}\hat{\mathscr{L}}_k^{\beta}(x)\hat{\mathscr{L}}_l^{\beta}(x)dx=\frac{1}{\beta}\delta_{kl}.
\end{equation}
We then define
\begin{equation}
\label{eq:basis_inf}
 \phi^{\infty}_j(z)=\hat{\mathscr L}^\beta_j(z-L), \ \ \ j=0,\dots,q
\end{equation}
for which the analog of \eqref{eq:normalize_phi} holds 
\begin{equation}
\int_{L}^{+\infty} \phi^{\infty}_{k}(z) \phi^{\infty}_{l}(z) \,dz = \frac1{\beta} \delta_{kl}.
\end{equation}
and we assume that
\begin{equation}
\label{eq:modalexplag}
c(z,t)\approx \sum_{j=0}^{q}  c^{(j)}_{\infty}(t) \phi^{\infty}_j(z) \ \ \ \ z\in K_{\infty}.
\end{equation}
For the resulting integrals, Gauss-Laguerre-Radau quadrature will be employed, see  \cite{benacchio:2013} for the definition. Approximation \eqref{eq:modalexplag} amounts to say that the restriction of the numerical approximation of $c $ to $K_{\infty} $ will be sought in the linear space $ V^q_{\infty} $ spanned by the functions defined in \eqref{eq:basis_inf}.\\Therefore, the global finite element space employed in the proposed extended DG discretization can be identified with 
$ V^{p,q}_{h} =  V^p_{h}\oplus V^q_{\infty}. $
For $v\in V^{p,q}_{h}$ we can then introduce the jump and average operators as (see, e.g., \cite{arnold:2002}) 
\begin{gather}
    \llbracket v(z)\rrbracket=v(z^-)-v(z^+), \qquad  
    \{v(z)\}=\frac{1}{2}(v(z^-)+v(z^+))
\end{gather}
and we remark that for $u,v \in V^{p,q}_{h} $ one has
\begin{equation}
    \llbracket uv \rrbracket= \left\{u\right\}\llbracket v\rrbracket+\left\{v\right\}\llbracket u\rrbracket.
\end{equation}
The extended DG discretization then involves integration of equation \eqref {eq:advdiffreacgen} against a test function $v\in V^{p,q}_{h}$, integrating by parts and imposing for $ m=1,\dots,N $ the appropriate continuity constraints at the interelement boundaries.
Setting
\begin{equation}
N_{\infty}=\{1,2,\dots,N,\infty\},     
\end{equation}
denoting by $\sigma >0$ the stabilization parameter for the DG approximation of the parabolic terms, and redefining the jump and average operators at $z=0$ so as to account for the boundary conditions, we obtain the following weak extended DG formulation of the problem:
\begin{align}
    \nonumber\sum_{m\in N_{\infty}}\int_{K_m}\frac{\partial c}{\partial t}vdz &=
    -\sum_{m=1}^{N}\hat f_{m+1/2}\left\llbracket v\left(z_{m+1/2}\right)\right\rrbracket  \\\nonumber
    &+\sum_{m\in N_{\infty}}\int_{K_m}f(c)v^\prime dz -\sum_{m\in N_{\infty}}\int_{K_m}\mu\frac{\partial c}{\partial z}v^\prime dz \\\nonumber
    &+\sum_{m=0}^N\left\{\mu\left(z_{m+1/2}\right)\frac{\partial c}{\partial z}\left(z_{m+1/2}\right)\right\}\left\llbracket v\left(z_{m+1/2}\right)\right\rrbracket \\\nonumber
    &+\sum_{m=0}^N\left\{\mu\left((z_{m+1/2}\right) v^\prime\left(z_{m+1/2}\right)\right\}\left\llbracket c\left(z_{m+1/2}\right)\right\rrbracket
   \\\nonumber
   &-\sum_{m=0}^N\frac{\sigma}{\Delta z_m}\left\llbracket c\left(z_{m+1/2}\right)\right\rrbracket\left\llbracket v\left(z_{m+1/2}\right)\right\rrbracket  \\\nonumber
    &+ \sum_{m\in N_{\infty}}\int_{K_m} svdz \\
    &+\mu(0,t) v^\prime(0)g_0(t)+f(g_0(t))v(0). 
\end{align}
At $z=L$,  the limit from the left of the approximate solution is computed as $\sum_{j=0}^{p}  c^{(j)}_{N}(t) \phi^N_j(L)$,
while the limit from the right  is computed as $\sum_{j=0}^{q}  c^{(j)}_{\infty}(t)$, according to the approximations \eqref{eq:modalexpdg}, \eqref{eq:modalexplag},
respectively. Among the many possible formulations for the parabolic terms, for definiteness we choose that corresponding to the Symmetric Interior Penalty Galerkin method (SIPG), see, e.g.,  \cite{arnold:1982, wheeler:1978} and the review in \cite{riviere:2008}. Furthermore, the Rusanov numerical flux is employed for the hyperbolic terms, so that
\begin{equation}\label{eq:rusanov}
  \hat f_{m+1/2} = \frac{1}{2}\left[f\left(c_{h,m+1/2}^+\right)+f\left(c_{h,m+1/2}^-\right)\right]
  -\frac{\Lambda_{m+1/2}}{2}\left(c_{h,m+1/2}^+-c_{h,m+1/2}^-\right),
\end{equation}
where the time dependency is omitted for simplicity,
\begin{equation}
\Lambda_{m+1/2}=\text{max}\left(\left\lvert \frac{df}{dc}\left(c_{h,m+1/2}^+\right)\right\rvert,\left\lvert \frac{df}{dc}\left(c_{h,m+1/2}^-\right)\right\rvert\right),
\end{equation}
and $c_{h,m+1/2}^+=c_h \left(z_{m+1/2}^+\right),\,c_{h,m+1/2}^-=c_h \left(z_{m+1/2}^-\right)$.

\noindent We can now define the bilinear form 
\begin{equation}
a:  V^{p,q}_{h}\times  V^{p,q}_{h}\times[0,+\infty)\rightarrow\mathbb{R} \rm \ \ as
\end{equation}
\begin{align}\label{bil_form}
\begin{split}
    a(w,v,t)=&\sum_{m\in N_{\infty}}\int_{K_m}\mu\frac{\partial w}{\partial z}v^\prime dz\\
    &-\sum_{m=0}^N\left\{\mu\left(z_{m+1/2},t\right)\frac{\partial w}{\partial z}\left(z_{m+1/2}\right)\right\}\left\llbracket v\left(z_{m+1/2}\right)\right\rrbracket\\&-\sum_{m=0}^N\left\{\mu\left(z_{m+1/2},t\right) v^\prime\left(z_{m+1/2}\right)\right\}\left\llbracket w\left(z_{m+1/2}\right)\right\rrbracket\\&+\sum_{m=0}^N\frac{\sigma}{\Delta z_m}\left\llbracket w\left(z_{m+1/2}\right)\right\rrbracket\left\llbracket v\left(z_{m+1/2}\right)\right\rrbracket,
    \end{split}
\end{align}
and the nonlinear function $b:  V^{p,q}_{h} \times V^{p,q}_{h}\rightarrow\mathbb{R}$ as
\begin{align}\label{nonlin_form}
\begin{split}
    b(w,v)&=\sum_{m=1}^{N}f\left(w\left(z_{m+1/2}\right)\right)\left\llbracket v\left(z_{m+1/2}\right)\right\rrbracket\\&-\sum_{m\in N_{\infty}}\int_{K_m}f(w)v^\prime dz-\sum_{m\in N_{\infty}}\int_{K_m} s(w)vdz.
    \end{split}
\end{align}
We also introduce $g:V^{p,q}_{h}\times V^{p,q}_{h}\to\mathbb{R}$, $h:V^{p,q}_{h}\times V^{p,q}_{h}\to\mathbb{R}$ as
\begin{align}
    g(w,v,t)&=-\mu(L,t) v^\prime(L)w(L)+\frac{\sigma}{\Delta z_N}v(L)w(L) \\
    h(w,v)&=-f(w(L))v(L)
\end{align}
and the linear operator $L: V^{p,q}_{h}\times[0,+\infty)\rightarrow\mathbb{R}$
\begin{equation}\label{lin_op}
    L(v,t) = \mu(0,t) v^\prime(0)g_0(t)+\frac{\sigma}{\Delta z_1}v(0)g_0(t)+f(g_0(t))v(0),
\end{equation}
which is related to the Dirichlet condition at the left endpoint $z=0$.  
The extended DG weak formulation can then be written more compactly as follows:\\

\noindent \textit{For all $t>0$, find $c_h(t)\in  V^{p,q}_{h}$ such that, $\forall v\in  V^{p,q}_{h}$,}
\begin{equation}\label{weak_form_DG}
    \int_0^{+\infty}\frac{\partial c_h}{\partial t}vdz=-a(c_h,v,t)-b(c_h,v)
    +L(v,t)+g(c_h,v,t)+h(c_h,v).    
\end{equation}
Approximating $c_h$ using \eqref{eq:modalexpdg} and \eqref{eq:modalexplag},
and taking  $v=\phi_j^m $ and $v=\phi_k^{\infty}$, 
one obtains a set of equations for the discrete degrees of freedom 
$c_m^{(j)} \ \  j=0,\dots,p, m=1,\dots,N $ $ c_{\infty}^{(k)} \  k=0,\dots,q.$
Collecting  these in two time-dependent vectors $ \mathbf{c}_{Dg}\in \mathbb{R}^{N(p+1)}  $  and $ \mathbf{c}_{Lg}\in \mathbb{R}^{q+1},$
one obtains the systems
\begin{align}\label{discr_form}
\begin{split}
    \frac{d\mathbf{c}_{Dg}}{dt}&=\mathbf{A}_{Dg}\mathbf{c}_{Dg}+\mathbf{A}_{Dg,Lg}\mathbf{c}_{Lg}
    +\mathbf{b}_{Dg}(\mathbf{c}_{Dg})+\mathbf{h}_{Dg}( \mathbf{c}_{Dg},\mathbf{c}_{Lg})  +\mathbf{g}_0(t) \end{split}\\\begin{split}
      \frac{d\mathbf{c}_{Lg}}{dt}&=\mathbf{A}_{Lg}\mathbf{c}_{Lg}+\mathbf{A}_{Lg,Dg}\mathbf{c}_{Dg}
    +\mathbf{b}_{Lg}(\mathbf{c}_{Lg}) +\mathbf{h}_{Lg} (\mathbf{c}_{Dg},\mathbf{c}_{Lg}).
    \end{split}
\end{align}
The time-dependent matrices $\mathbf{A}_{Dg}$ and $\mathbf{A}_{Lg} $ result from the discretization of the diffusion operator in the interior of the $[0,L] $ and $[L,+\infty) $ subdomains, respectively. The coupling matrices $\mathbf{A}_{Dg,Lg}$ and
$\mathbf{A}_{Lg,Dg}$ result from the discretization
of the diffusion operator involving discrete degrees of freedom of both subdomains. The nonlinear functions
$\mathbf{b}_{Dg}$ and $\mathbf{b}_{Lg} $ result from the discretization
of the hyperbolic part and source terms in the interior of the $[0,L] $ and $[L,+\infty) $ subdomains, respectively. 
The term $\mathbf{g}_0(t) $ is associated with boundary conditions at $z=0$, while the coupling nonlinear functions $\mathbf{h}_{Dg}( \mathbf{c}_{Dg},\mathbf{c}_{Lg}) $ and $ \mathbf{h}_{Lg}( \mathbf{c}_{Dg},\mathbf{c}_{Lg}) $ result from the discretization of the hyperbolic part involving discrete degrees of freedom of both subdomains.\\
Next, we define the global unknown vector as
\begin{equation}
    \mathbf{c}(t)=\left(\mathbf{c}_{Dg}(t), \mathbf{c}_{Lg}(t)\right)^T\in\mathbb{R}^{N(p+1)+q+1},
\end{equation}
and the global vectors
\begin{gather}
    \mathbf{b}\left(\mathbf{c}(t)\right)=\left(\mathbf{b}_{Dg}\left(\mathbf{c}_{Dg}(t)\right), \mathbf{b}_{Lg}\left(\mathbf{c}_{Lg}(t)\right)\right)^T \in \mathbb{R}^{N(p+1)+q+1}\\
    \mathbf{h}(\mathbf{c})=\left(\mathbf{h}_{Dg}\left( \mathbf{c}_{Dg},\mathbf{c}_{Lg}\right), \mathbf{h}_{Lg}\left( \mathbf{c}_{Dg},\mathbf{c}_{Lg}\right)\right)^T \in \mathbb{R}^{N(p+1)+q+1}\\
    \mathbf{g}(t)=\left(\mathbf{g}_{Dg}(t),0,\dots,0\right)^T \in \mathbb{R}^{N(p+1)+q+1}.
\end{gather}
Defining the global extended DG matrix
\begin{equation}\label{eq:globalmatrix}
    \mathbf{A}(t)=\begin{pmatrix}
                   \mathbf{A}_{Dg}(t) & \mathbf{A}_{Dg,Lg}(t)\\\
                   \mathbf{A}_{Lg,Dg}(t) & \mathbf{A}_{Lg}(t)
               \end{pmatrix}\in\mathbb{R}^{(N(p+1)+q+1)\times (N(p+1)+q+1)},
\end{equation}
the extended DG semi-discrete formulation reads
\begin{equation}\label{eq:semidiscrglobal}
    \frac{d\mathbf{c}(t)}{dt}=\mathbf{A}(t)\mathbf{c}(t)+\mathbf{b}(\mathbf{c}(t))+\mathbf{h}(\mathbf{c}(t))+\mathbf{g}(t).
\end{equation}
The matrix $\mathbf{A}$ is the discretization of the diffusion term, the vector $\mathbf{b}$ is the discretization of the non-linear advective part and the optional source-reaction term, the vector $\mathbf{h}$ contains the flux exchange at the interface $z=L$ by means of the application of the Rusanov flux to the flux function $f$, and the vector $\mathbf{g}$ encodes the Dirichlet condition at the left endpoint $c(0)=a(t)$. 
We remark that, because of the vectors $\mathbf{b}$ and $\mathbf{h}$, problem \eqref{eq:semidiscrglobal} is non-linear. However, if the functions $f$ and $s$ in \eqref{eq:advdiffreacgen} are linear, then $\mathbf{b}(\mathbf{c}(t))$ and $\mathbf{h}(\mathbf{c}(t))$ can be written as the product between a matrix and the unknown vector $\mathbf{c}(t)$. In this case, \eqref{eq:semidiscrglobal} is a linear system of equations.

\noindent The semi-discrete extended DG formulation \eqref{eq:semidiscrglobal} can then be  discretized in time by any standard method for the numerical solution of ODE systems. In this paper, we use the Crank-Nicolson method for the linear test problems considered in Section \ref{sec:experiments}. For the non-linear problems, a second order implicit-explicit (IMEX) method is used, that is described, e.g., in \cite{bonaventura:2017, giraldo:2013}. \revision{As the terms associated with the diffusion process can entail rather restrictive stability constraints on the time step size if discretized explicitly, they are discretized implicitly, while the terms associated with the hyperbolic conservation law are treated explicitly.}
 
\section{Stability analysis}
\label{sec:analysis} \indent  

\noindent In order to study the numerical stability of the global semi-discrete extended DG formulation \eqref{eq:semidiscrglobal}, we consider the special case of the linear, constant coefficient, advection-diffusion equation:
\begin{equation}\label{eq:advdifflin}
    \frac{\partial c}{\partial t}+u\frac{\partial c}{\partial z}=\mu\frac{\partial^2 c}{\partial z^2}+s(c,z,t)
\end{equation}
i.e. equation \eqref{eq:advdiffreacgen} with $f(c)=uc$, $\mu(z,t)\equiv\mu$, $\mu,u\in\mathbb{R}^+$, and in the homogeneous case the source term $s\equiv0$.
\revision{ We also assume for definiteness $u>0.$ The application of the extended DG scheme results in the semi-discretization
\begin{equation}
    \frac{d\mathbf{c}(t)}{dt}=\mathbf{A}\mathbf{c}(t) + \mathbf{f}(t)
    \label{eq:semi_disc}
\end{equation}
We describe the structure of the matrix $\mathbf{A}$ below.
Note that, in this section, the matrix $\mathbf{A}$ includes both the advective and diffusive terms.

\begin{table}[h]\centering
\begin{tabular}{|clllll|ll|llll|}
\hline
\multicolumn{8}{|c|}{\multirow{8}{*}{$\mathbf{A}_{Dg}$}} & & & &  \\
\multicolumn{8}{|c|}{} & & & &  \\
\multicolumn{8}{|c|}{} & & & &  \\
\multicolumn{8}{|c|}{} & & & &  \\
\multicolumn{8}{|c|}{} & & & &  \\ \cline{9-12} 
\multicolumn{8}{|c|}{} & & \multicolumn{2}{l}{\multirow{3}{*}{$\mathbf{A}_{Dg,Lg}$}} &  \\
\multicolumn{8}{|c|}{} & & \multicolumn{2}{l}{} &  \\
\multicolumn{8}{|c|}{} & & \multicolumn{2}{l}{} &  \\ \hline
\multicolumn{1}{|l}{}  & & & & & & & & & & &  \\
\multicolumn{1}{|l}{}  & & & & & & \multicolumn{2}{l|}{\multirow{3}{*}{$\mathbf{A}_{Lg,Dg}$}} &  & \multicolumn{2}{c}{\multirow{3}{*}{$\mathbf{A}_{Lg}$}}    &  \\
\multicolumn{1}{|l}{}  & & & & & & \multicolumn{2}{l|}{} & & \multicolumn{2}{c}{} &  \\
\multicolumn{1}{|l}{} & & & & & & \multicolumn{2}{l|}{} & & \multicolumn{2}{c}{} &  \\
\multicolumn{1}{|l}{} & & & & & & & & & & &  \\ \hline
\end{tabular}
\end{table}

\noindent The matrix $\mathbf{A}$ consists of four blocks, $\mathbf{A}_{Dg}$,$\mathbf{A}_{Lg}$, $\mathbf{A}_{Dg,Lg}$ and $\mathbf{A}_{Lg,Dg}$.
For the sake of simplicity, we assume that   $\Delta z_m=\Delta z$ for all $m=1,\dots,N$; in this case, the boundary quantities are independent of the element index $m$ and can be denoted as
\begin{align}
    \phi_i^L&=\phi_i^m(z_{m-\frac{1}{2}}^+) \quad &\phi_i^R&=\phi_i^m(z_{m+\frac{1}{2}}^-)  \\ (\phi_i^\prime)^L&=(\phi_i^m)^\prime(z_{m-\frac{1}{2}}^+) \quad &(\phi_i^\prime)^R&=(\phi_i^m)^\prime(z_{m+\frac{1}{2}}^-)
\end{align}

\noindent The DG discretization on $[0,L]$ is described by
the block tridiagonal matrix

\begin{equation}
    \mathbf{A}_{Dg}=-\frac{1}{\Delta z}\begin{bmatrix} B_1 & D & & & \\ R_2+E & B_2 & D & & \\ & \ddots & \ddots & \ddots & & \\ & & R_2+E & B_2 & D \\ & & & R_2+E & B_2\end{bmatrix}\in\mathbb{R}^{N(p+1)\times N(p+1)}
\end{equation}
 each block being of dimension $(p+1)\times (p+1)$, where $B_1=A+F+C+DG_{adv}+R_1$ and $B_2=A+B+C+DG_{adv}+R_1$, $i,j=0,\dots,p$,

\begin{gather}
    A_{i+1,j+1}=\mu\int_{K_m}\phi_i^\prime(z)\phi_j^\prime(z)\,dz \\\nonumber \\ 
    B_{i+1,j+1}=\frac{\mu}{2}(\phi_j^\prime)^L\phi_i^L - \frac{\epsilon\mu}{2}\phi_j^L(\phi_i^\prime)^L+\frac{\sigma}{\Delta z}\phi_j^L\phi_i^L \\ \nonumber\\
    C_{i+1,j+1}=-\frac{\mu}{2}(\phi_j^\prime)^R\phi_i^R + \frac{\epsilon\mu}{2}\phi_j^R(\phi_i^\prime)^R+\frac{\sigma}{\Delta z}\phi_j^R\phi_i^R \\\nonumber \\
    D_{i+1,j+1}=-\frac{\mu}{2}(\phi_j^\prime)^L\phi_i^R - \frac{\epsilon\mu}{2}\phi_j^L(\phi_i^\prime)^R-\frac{\sigma}{\Delta z}\phi_j^L\phi_i^R \\\nonumber \\
    E_{i+1,j+1}=\frac{\mu}{2}(\phi_j^\prime)^R\phi_i^L + \frac{\epsilon\mu}{2}\phi_j^R(\phi_i^\prime)^L-\frac{\sigma}{\Delta z}\phi_j^R\phi_i^L \\\nonumber \\
    F_{i+1,j+1}=-\mu(\phi_j^\prime)^R\phi_i^R + \epsilon\mu\phi_j^R(\phi_i^\prime)^R+\frac{\sigma}{\Delta z}\phi_j^R\phi_i^R \\\nonumber \\
    (DG_{adv})_{i+1,j+1}=-u\int_{K_m}\phi_j(z)\phi^\prime_i(z)\,dz \\\nonumber\\
    (R_1)_{i+1,j+1}=u\phi_j^R\phi_i^R \\\nonumber \\
    (R_2)_{i+1,j+1}=-u\phi_j^R\phi_i^L
\end{gather}
 The semi-infinite Laguerre discretization is instead described by
\begin{equation}
    \mathbf{A}_{Lg}=-\mu\beta^2L_fL_f^T + u\beta L_f + LAG_{Dg} \in\mathbb{R}^{(q+1)\times(q+1)}
\end{equation}
where $L_f\in\mathbb{R}^{(q+1)\times (q+1)}$ has $-1/2$ on the diagonal, $-1$ in the lower triangular part and $0$ in the upper triangular part, while
\begin{equation}
    (LAG_{Dg})_{i+1,j+1}=\beta\left[\frac{\mu\beta}{2}\left(j+\frac{1}{2}\right)-\frac{\mu\epsilon\beta}{2}\left(i+\frac{1}{2}\right)-\frac{\sigma}{\Delta z}\right] \quad i,j=0,\dots,q.
\end{equation}
Unlike $\mathbf{A}_{Dg}$, $\mathbf{A}_{Lg}$ is a full matrix.

\noindent The coupling between the finite and semi-infinite subdomain in the extended DG scheme is represented by the matrices $\mathbf{A}_{Dg,Lg}\in\mathbb{R}^{(p+1)\times(q+1)}$ and $\mathbf{A}_{Lg,Dg}\in\mathbb{R}^{(q+1)\times(p+1)}$, defined as
\begin{align}
    (\mathbf{A}_{Dg,Lg})_{i+1,j+1}&=-\frac{\mu\beta}{2\Delta z}\phi_i^R\left(j+\frac{1}{2}\right)+\frac{\mu\epsilon}{2\Delta z}(\phi_i^\prime)^R+\frac{\sigma}{\Delta z^2}\phi_i^R \\\nonumber &i=0,\dots,p \quad j=0,\dots,q
    \end{align}
    \begin{align}
    (\mathbf{A}_{Lg,Dg})_{i+1,j+1}&=-\frac{\mu\beta}{2}(\phi_j^\prime)^R+\frac{\mu\epsilon\beta^2}{2}\phi_j^R\left(i+\frac{1}{2}\right)+\frac{\sigma\beta}{\Delta z}\phi_j^R+u\beta\phi_j^R \\\nonumber& i=0,\dots,q \quad j=0,\dots,p
\end{align}

\noindent Proving that $\mathbf{A}$ has eigenvalues with negative real part in the most general case is not immediate, but this can be achieved in a rather straightforward way in the purely advective, inviscid case.

\begin{theorem}
If $\mu=0, $ $\mathbf{A}$ has eigenvalues with negative real part for all values of $N,q,p.$
\end{theorem}
\begin{proof}
In the inviscid case $\mu=0, $ one has $\mathbf{A}_{Dg,Lg} = \mathbf{0},$  so that $\mathbf{A}$ is a block lower triangular matrix. Therefore, its eigenvalues coincide with those of the blocks that include the main diagonal.
One of these blocks is $\mathbf{A}_{Lg}, $
whose eigenvalues are all equal to $-u\beta/2.$ The other block is itself a block lower triangular matrix, whose eigenvalues are given by the eigenvalues of $B_2$ taken with a multiplicity equal to the number of elements $N.$ These eigenvalues  can be computed directly and shown to have negative real part independently of $N,p. $ For example, in the $p=1 $ case they are given by
\begin{equation}
\lambda_{\pm}=
\left(-2 \pm i\sqrt{2}\right)\frac{u}{\Delta z}.    
\end{equation}
\qed
\end{proof}

\noindent In order to provide an empirical check of the stability of this formulation also in the diffusive case, we compute the spectrum of the extended DG matrix $\mathbf{A}$ as a function of the P\'eclet number $Pe=u\mathcal{L}/\mu$, where $\mathcal{L}$ is a reference length scale, for fixed values of $N $ and $q.$ More specifically, we set $\mathcal{L}=1$, $\beta=1$ for the scaling Laguerre parameter,  $\mu=1,\, u=Pe\mu$. The polynomial degree  $p=2$ is used in the DG discretization in the finite subdomain. Results of the analysis are reported in Table \ref{Tab:pe_fullyDG}. For all values of $Pe$, all eigenvalues have negative real part, giving empirical corroboration to the stability of the extended DG scheme. While only results for $\epsilon=1$ are shown in Table \ref{Tab:pe_fullyDG}, qualitatively equivalent figures are obtained for $\epsilon=0,-1$ and different number of elements in the finite subdomain (not shown).
\begin{table}[htbp]\centering\footnotesize
\caption{Maximum real part of the eigenvalues of $A$ and relative $L^2$ and $L^\infty$ errors as a function of $Pe$ of the extended DG scheme with respect to the exact solution, linear non-homogeneous advection-diffusion equation. $N=100$, $q=180$, $\Delta t=0.05\,\textrm{s}$, $u=1$, $\mu=u/Pe$, $\sigma=200$, $\epsilon=1$. For the advective case $Pe=\infty$, $\sigma=0$.}\label{Tab:pe_fullyDG}
 \centering
\begin{tabular}{cccc}\toprule\midrule
 $Pe$             & $Re(\lambda)_{max}$ & $\mathcal{E}_2^{\textrm{rel}}$  & $\mathcal{E}_\infty^{\textrm{rel}}$\\\midrule
 0.001 & -1.90E-02 & 1.58E-03 & 6.97E-04 \\
 10    & -2.13E-02 & 3.39E-03 & 2.83E-03 \\
 100   & -2.41E-02 & 2.79E-03 & 2.74E-03 \\
 500   & -2.57E-02 & 2.68E-03 & 2.69E-03 \\
 1000  & -2.61E-02 & 2.66E-03 & 2.69E-03 \\
 10000 & -1.66E-02 & 2.65E-03 & 2.68E-03 \\
 100000 & -1.66E-03 & 2.64E-03 & 2.68E-03 \\
 1000000 & -1.73E-04 & 2.56E-03 & 2.66E-03 \\
 $\infty$ & -5.00E-01 & 2.65E-03 & 2.65E-03 \\
 \bottomrule\end{tabular}
 \end{table}}

 \section{Numerical experiments}
\label{sec:experiments} \indent  
We present here the results of several numerical tests with the extended DG approach described in the previous sections. First, a number of validation tests are carried out, considering both linear and non-linear model problems. The tests assess the accuracy of the extended DG scheme by comparing it with a stand-alone, single-domain reference DG discretization on a wider domain. Errors are computed using a relatively large number of modes in the semi-infinite subdomain of the extended DG scheme, and the wider domain for the stand-alone reference run also covers part of the semi-infinite subdomain.
Next, we add a damping reaction term in the semi-infinite subdomain to simulate an absorbing layer. We show that the extended DG scheme efficiently damps signals leaving the finite subdomain with negligible reflections into the finite region \revision{as compared with a damped single-domain DG scheme using both uniform and non-uniform grids}.

The experiments consider the linear advection-diffusion  equation with constant coefficients \eqref{eq:advdifflin}, both in the non-homogeneous case ($s\neq 0$) and the homogeneous case ($s=0$), and the nonlinear, homogeneous viscous Burgers' equation with constant viscosity, i.e., equation \eqref{eq:advdiffreacgen} with $f(c)=c^2/2$, $\mu(z,t)\equiv\mu\in\mathbb{R}^+$, and $s=0$.   
Errors are computed on the finite region $[0,L]$ using a suitable Gaussian quadrature rule on the sub-intervals $K_m$, whose width is $\Delta z=L/N$ for all $m=1,\dots,M$. In particular, we introduce the discrete norms
\begin{align}
    \lVert c_h\rVert_{L^2}&=\sqrt{\sum_{m=1}^N\frac{\Delta z}{2}\sum_{k=1}^{ng}\left[c_h\left(\frac{\Delta z}{2}x_k+z_m\right)\right]^2w_k} \\
    \lVert c_h\rVert_{L^\infty}&=\max_{m=1,\dots,N}\max_{k=1,\dots,ng}\left\lvert c_h\left(\frac{\Delta z}{2}x_k+z_m\right)\right\rvert,
\end{align}
where $\{x_k\}_{k=1}^{ng}$ and $\{w_k\}_{k=1}^{ng}$ are the Gaussian nodes and weights on the reference interval $[-1,1]$, with $ng$ the number of quadrature points. Absolute errors with respect to a reference solution are defined as
\begin{equation}
    \mathcal{E}_r=\lVert c_h-c_{ref}\rVert_{L^r},\quad r\in{2,\infty}
\end{equation}
where $c_h$ and $c_{ref}$ are the numerical and the reference solution, respectively; the latter may be either the exact solution or a single-domain DG discretization. In some tests we will be also interested in relative errors with respect to the reference solution $c_{ref}$ defined as

\begin{equation}
\mathcal{E}_r^{rel}=\frac{\lVert c_h-c_{ref}\rVert_{L^r}}{\lVert c_{ref}\rVert_{L^r}},\quad r\in{2,\infty}.
\end{equation}

\newpage

\subsection{Validation of the extended DG scheme coupling strategy}
We start by testing the proposed method for the linear advection-diffusion equation with constant coefficients \eqref{eq:advdifflin}. We first consider the non-homogeneous case by setting $s\neq 0$,  assuming an exact solution $c(z,t)= ze^{-z}\text{sin}^2(z-t)$, and computing the right-hand side analytically.
 We run the scheme for a variable number of modes $q$ in the semi-infinite subdomain, also setting $N=100$, $L=2\,\textrm{m}$, $p=2$, $\mu=1\,\textrm{m}^2/\textrm{s}$, $u=2Pe\mu$, and a final time $T=10\,\textrm{s}$ with $n=200$ time steps. In order to minimize errors at the interface, the value of $\beta$ is chosen in such a way that the distance between the first and the second node in the semi-infinite subdomain matches the grid spacing $\Delta z$ in the finite subdomain. Relative $L^2$ and $L^\infty$ errors at time $T$ with respect to the exact solution due to the use of different sets of basis functions in the finite and semi-infinite subdomains are below $5\times10^{-6}$ using at least 20 modes, below $3\times10^{-3}$ using 10 modes, and a few percent using 5 modes, thus displaying spectral convergence in space  (Table \ref{Tab:M_fullyDG}\revision{, scaling parameter $\beta$ chosen for matching grid spacing at the finite/semi-infinite interface).}  

\begin{table}[htbp]\footnotesize
\caption{Relative $L^2$ (\revision{$\mathcal{E}_2^{rel}$}) and $L^\infty$ (\revision{$\mathcal{E}_\infty^{rel}$}) errors of the extended DG scheme with respect to the exact solution, linear non-homogeneous advection-diffusion equation, \revision{several values for the number of modes $q$ and scaling parameter $\beta$}. $\Delta t = 0.005\,\textrm{s}$, $\Delta z=0.02\,\textrm{m}$, $C=0.25$.}\label{Tab:M_fullyDG}
\centering
\begin{tabular}{cccc}\toprule\midrule
$q$ & $\beta$             & $\mathcal{E}_2^{\textrm{rel}}$ & $\mathcal{E}_\infty^{\textrm{rel}}$\\\midrule
5 & 30 & 5.39E-02 & 7.93E-02 \\
10 & 16 & 2.39E-03 & 3.24E-03 \\
20 &  8 & 3.35E-06 & 2.99E-06 \\
40 &  4 & 3.35E-06 & 2.99E-06 \\
80 &  2 & 3.35E-06 & 2.99E-06 \\
\bottomrule\end{tabular}
\end{table}
\revision{\noindent To further inspect the properties of the extended DG scheme, we evaluate the relative $L^2$, and $L^\infty$ errors for varying $\beta$ values (Figure \ref{fig:beta_study}) and $q=40$ (left panel), $q=5$ (right panel). The markers show the error value for the choice of $\beta$ corresponding to the matching of the spacing at the interface,  $\hat{z}_2-\hat{z}_1\approx\Delta z$, obtained with $\beta=30$ and $q=5$ modes, and $\beta=4$ and $q=40$ modes. For $q=40$,  errors are constant for $\beta\leq20$ and lowest around $\beta=25$. By contrast, for $q=5$, errors are lowest around $\beta=5$, so $\beta=4$ is a particularly good choice. These tests show that an optimal choice of $\beta$ exists, but it is not necessarily the one determined by the matching condition $\hat{z}_2-\hat{z}_1\approx\Delta z$ at the interface between the finite and semi-infinite subdomains.

\begin{figure}[h]
    \centering
    \includegraphics[width=1.05\textwidth]{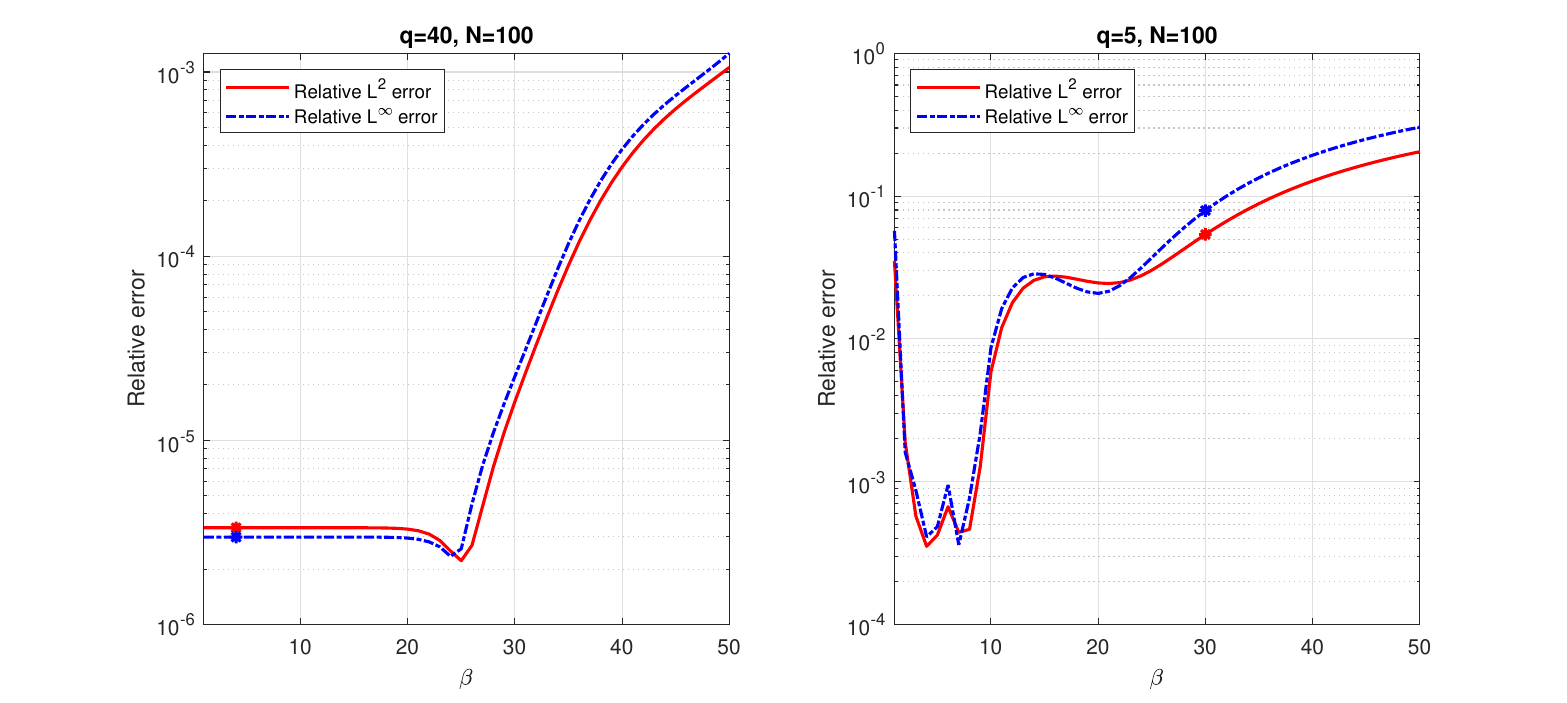}
    \caption{Relative $L^2$ and $L^\infty$ errors as a function of the scaling parameter $\beta$ in the extended DG scheme, linear non-homogeneous advection-diffusion equation. Stars: $\beta=4$ (left), $\beta=30$ (right).}
    \label{fig:beta_study}
\end{figure}
}

\begin{figure}
\centering
\includegraphics[width=.5\textwidth]{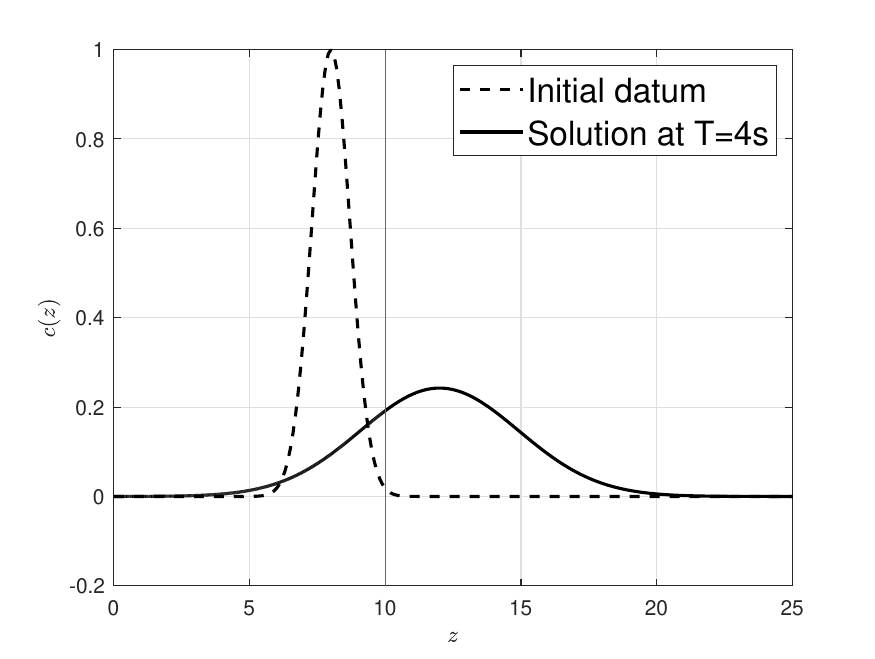}
\caption{Numerical solution of the linear homogeneous advection-diffusion equation. Dashed line: initial datum. Solid line: numerical solution with the extended DG scheme at $T=4\,\textrm{s}$. $z_c=8\,\textrm{m}$, $\sigma_c=1\,\textrm{m}$, $q=40$, $\beta=4$, $\Delta z=0.02\,\textrm{m}$, $\Delta t=0.02\,\textrm{s}$, $\mu=1\,\textrm{m}^2/\textrm{s}$, $u=1\,\textrm{m}/\textrm{s}$.}
\label{fig:hom_coupling}
\end{figure}

\noindent Next, we validate the extended DG approach in the homogeneous case. We consider a Gaussian initial datum:
\begin{equation}\label{eq:Gaussian}
    c_0(z)=\text{exp}\left[-\left(\frac{z-z_c}{\sigma_c}\right)^2\right].
\end{equation}
The interface is located at $L=10\,\textrm{m}$ and the initial hump is placed inside the bounded interval $[0,L]$ by choosing $z_c=8\,\textrm{m}$. The velocity is $u=1\,\textrm{m}/\textrm{s}$ and the final time is $T=4\,\textrm{s}$, so that the peak of the Gaussian crosses the interface, and the other parameters are $\Delta t=0.02\,\textrm{s}$, $q=40$ modes in the semi-infinite region, $N=500$ sub-intervals for the DG scheme, so that $\Delta z=0.02\,\textrm{m}$. As the model evolves, the initial hump expands and its amplitude decreases because of diffusion (Figure \ref{fig:hom_coupling}).
Relative errors in the finite subdomain $[0,10\,\textrm{m}]$ are computed for the extended DG scheme with respect to a single-domain DG solution run on $[0,50\,\textrm{m}]$. For~$q=10$ modes in the semi-infinite subdomain, relative errors are below a few percent, while for $q=40$ they lower to around $10^{-10}$ (Table \ref{tab:test2relerrs10_fullyDG}). 

\begin{table}[htbp]\footnotesize\caption{Relative $L^2$ (\revision{$\mathcal{E}_2^{rel}$}) and $L^\infty$ (\revision{$\mathcal{E}_\infty^{rel}$}) in $[0,10\,\textrm{m}]$ of the extended DG discretization with respect to a single-domain DG discretization, linear homogeneous advection-diffusion equation.  $z_c=8\,\textrm{m}$, $\Delta z = 0.02\,\textrm{m}$, $\Delta t = 0.02\,\textrm{s}$, $C=1$. Matching grid spacing at the interface is obtained for the chosen values of the scaling parameter $\beta$. }\label{tab:test2relerrs10_fullyDG}
\centering
\begin{tabular}{ccccc}\toprule\midrule
$q$ & $\beta$ & $\sigma_c$             & $\mathcal{E}_2^{\textrm{rel}}$ & $\mathcal{E}_\infty^{\textrm{rel}}$\\\midrule
\multirow{3}{*}{10} & \multirow{3}{*}{16} & 1 & 1.90E-02 & 3.80E-02 \\
 & & 2 & 1.98E-02 & 4.10E-02 \\
 & & 0.5 & 1.87E-02 & 3.71E-02 \\
\cmidrule(l{-.25pt}r{-.25pt}){1-5}
\multirow{3}{*}{40} & \multirow{3}{*}{4} & 1 & 3.51E-09 & 5.44E-08 \\
  &    & 2 & 4.30E-10 & 7.07E-09 \\
  &   & 0.5 & 6.51E-12 & 8.46E-11 \\
\bottomrule\end{tabular}
\end{table}

\begin{figure}[htbp]
\centering
\includegraphics[width=.48\textwidth]{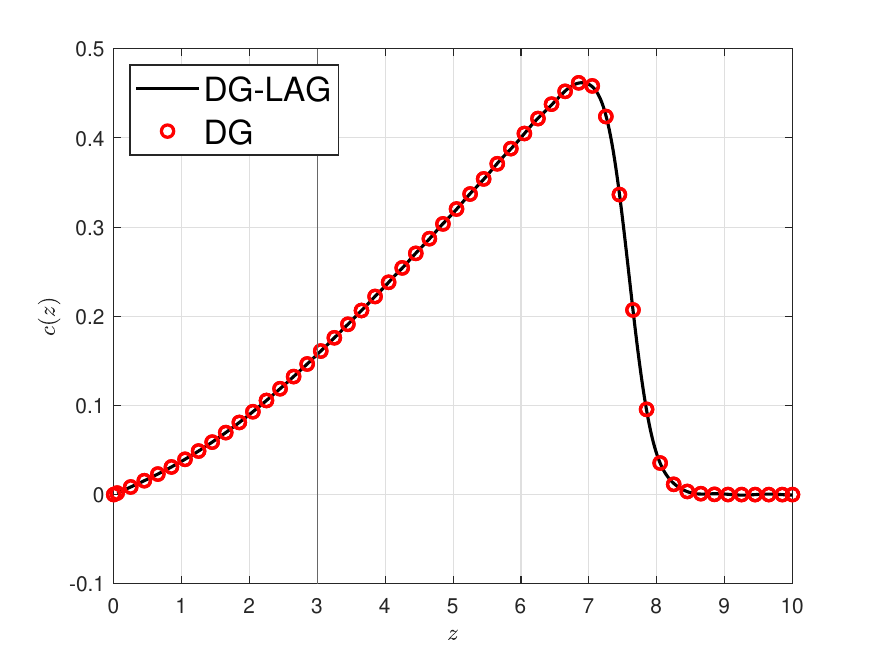}
\caption{Numerical solution at $T=1\,\textrm{s}$ of the viscous Burgers' equation with Gaussian initial datum $q_0(z)=\text{exp}\left(-(z-3)^2\right)$, $N=30$, $q=60$, $\Delta t=10^{-6}\,\textrm{s}$. Black line: extended DG scheme (DG-LAG). Red circles: Single-domain DG discretization.}
    \label{fig:burgerscoupled}
\end{figure}

\noindent In a final validation test, we consider the case of the viscous Burgers' equation \revision{with two different initial data. First, we consider t}he Gaussian profile $c_0(z)=\text{exp}\left(-(z-3)^2\right)$, with the interface placed at $z=3\,\textrm{m}$. As time $t$ evolves, the profile moves rightwards increasing its steepness -- using the viscosity value $\mu=0.05\,\textrm{m}^2/\textrm{s}$ no oscillations are observed until the final time $T=10\,\textrm{s}$. The extended DG scheme compares well with a reference solution computed by a single-domain DG discretization on $[0,10\,\textrm{m}]$ (Figure \ref{fig:burgerscoupled}).  

\noindent Since the closed form of the solution is not available, we compute the errors with respect to a stand-alone DG discretization on a larger domain with the same spacing $\Delta z$. The finite/semi-infinite interface in the extended DG scheme is placed at $L=3\,\textrm{m}$ and the model is run until $T=10\,\textrm{s}$, with $\Delta t=10^{-2}\,\textrm{s}$. The cases of $N=15$ and $N=30$ subintervals in the finite subdomain are considered, varying the number of modes in the semi-infinite subdomain, and, accordingly, the scaling parameter $\beta$ so that the distance between the first two nodes matches the grid spacing in $[0,L]$. The stand-alone single-domain DG reference solution for error computation is computed on the interval $[0,10\textrm{m}]$. A small number of Laguerre modes are found to suffice to keep the coupling errors in the bounded subdomain $[0,L]$ below a few percent (Table \ref{Tab:burgercoup_fullyDG}). 

\begin{table}[htbp]\footnotesize\caption{Relative $L^2$ (\revision{$\mathcal{E}_2^{rel}$}) and $L^\infty$ (\revision{$\mathcal{E}_\infty^{rel}$}) errors in $[0,3\,\textrm{m}]$ of the extended DG discretization with respect to a single-domain DG discretization, viscous Burgers' equation\revision{, for two different values of the grid spacing $\Delta z=L/N$ in the finite subdomain}. $T=10\,\textrm{s}$, $\Delta t=10^{-2}\,\textrm{s}$.  Matching grid spacing at the interface is obtained for the chosen values of the scaling parameter $\beta$.}\label{Tab:burgercoup_fullyDG}
\centering
\begin{tabular}{ccccc}\toprule\midrule
$N$ & $q$ & $\beta$             & $\mathcal{E}_2^{\textrm{rel}}$ & $\mathcal{E}_\infty^{\textrm{rel}}$\\\midrule
\multirow{4}{*}{15} & 10 & 1.6 & 2.10E-02 & 5.75E-02 \\
                    & 20 & 0.85 & 2.61E-02 & 6.70E-02 \\
                    & 40 & 0.45 & 2.72E-02 & 6.56E-02 \\
                    & 80 & 0.23 & 2.69E-02 & 6.35E-02 \\
\midrule
\multirow{4}{*}{30} & 10 & 3.6 & 6.21E-04 & 8.50E-04 \\
                    & 30 & 1.2 & 6.06E-04 & 1.29E-03 \\
                    & 60 & 0.6 & 6.77E-04 & 1.38E-03 \\
                    & 100 & 0.36 & 7.08E-04 & 1.36E-03 \\
\bottomrule\end{tabular}
\end{table}

\begin{figure}[h]
    \centering
    \includegraphics[width=0.5\textwidth]{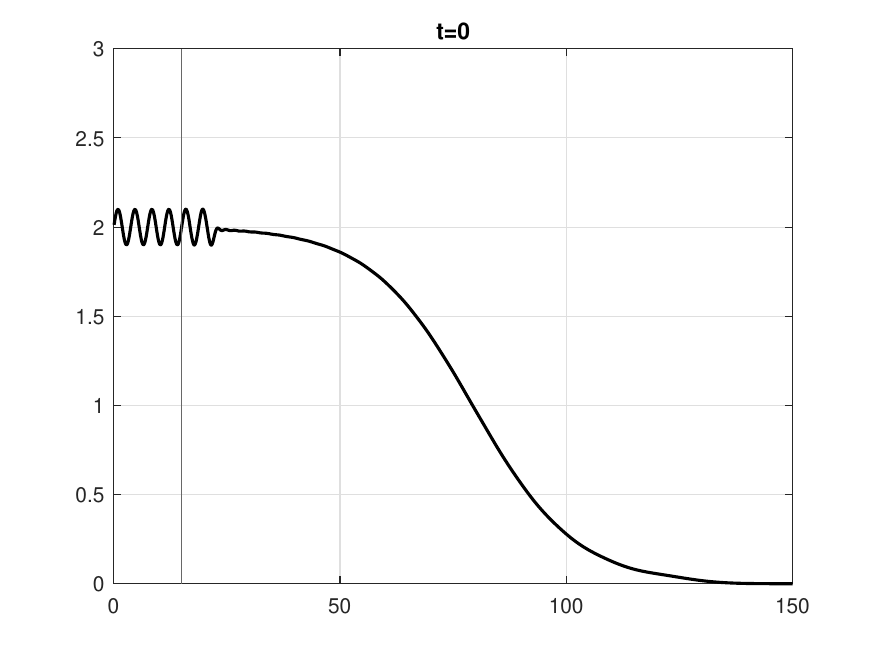}
    \caption{Initial profile for the sinusoidal validation test, $N=100$, $q=60$, $\beta=2$}
    \label{fig:initial_profile_burgers_wave}
\end{figure}

\revision{\noindent A further validation test with Burgers' equation aims to assess the ability of the extended DG scheme in simulating wave dynamics. The initial datum for this test is
\begin{equation}
\displaystyle
u_0(z)=\begin{cases}
C + A\sin{\left(\frac{k\pi}{L_0}z\right)} & z\leq L_0 \\
C\left(1-\frac{1}{1+\exp\left(\frac{\alpha\Bar{L}-(x-L_0)}{s}\right)}\right) & z>L_0
\end{cases}
\end{equation}
where $C=2\,m$, $A=0.1\,m$, $k=8$, $L=15\,m$, $\alpha=0.1$, $\Bar{L}=568.1231\,m$, $s=\Bar{L}/50$, $L_0=22.5\,m$ (Figure \ref{fig:initial_profile_burgers_wave})}.

\begin{figure}[h]
    \centering
    \includegraphics[width=\textwidth]{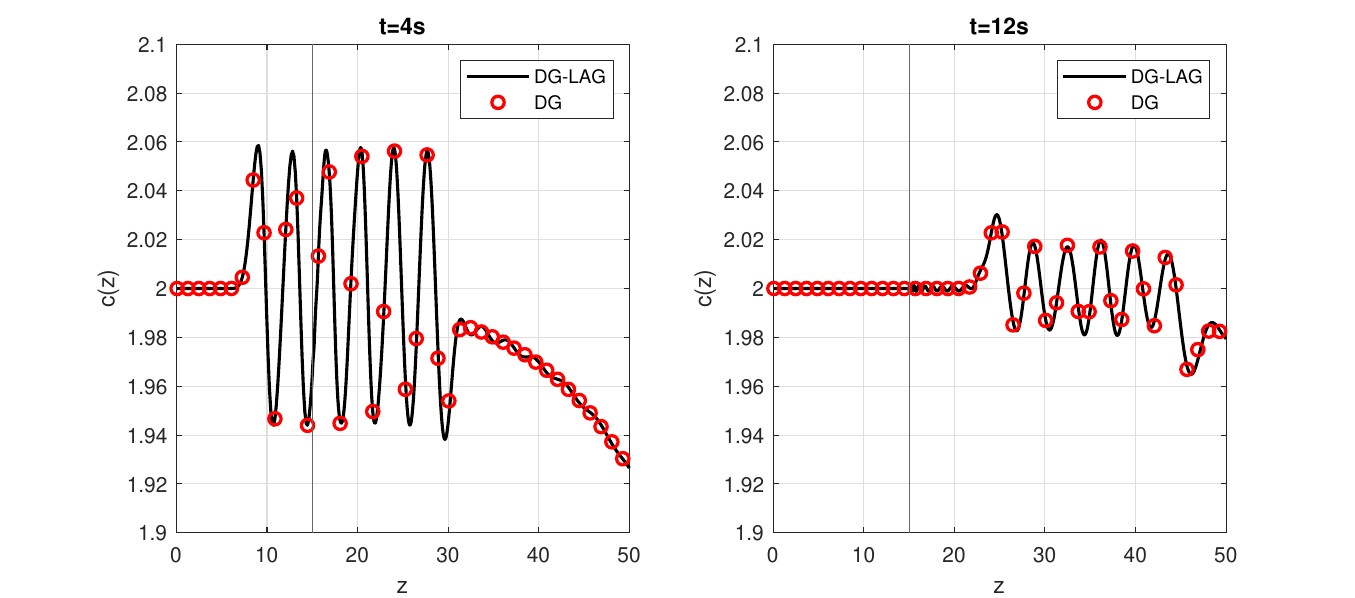}
    \caption{Sinusoidal initial data: extended DG solution (DG-LAG, solid line) and single-domain DG solution (DG, red dots) of viscous Burgers' equation at times $t=4,\,8,\,12\,\textrm{s}$, simulated with $\Delta t=10^{-2}\,\textrm{s}$.}
    \label{fig:burgers_wave}
\end{figure}

\revision{\noindent The extended DG scheme as well as a reference standalone DG scheme on the larger domain $[0,10L]=[0,150\,m]$ are run with $\Delta t=0.01\,s$, $N=100$, $q=60$, $\beta=2$.  
At final time $T=12\,s$, all the crests have crossed the interface placed at $L=15\,m$, and the extended DG scheme (solid black line) yields solutions close to the reference standalone DG scheme (red dots), see  
Figure \ref{fig:burgers_wave}. This is quantified by relative errors in the finite subdomain at final time of the extended DG scheme with respect to the reference solution (Table \ref{tab:burgers_wave_errors_const_beta} for constant scaling parameter $\beta$ and Table \ref{tab:burgers_wave_errors_matching_beta} for $\beta$ chosen to match the grid spacing at the interface). It is to be noted that the scaling parameter can be tuned to obtain lower error values. As previously observed, the value of $\beta$ corresponding to a matching grid spacing at the interface is not the optimal choice.}

\revision{\begin{table}[htbp]\footnotesize\caption{Relative $L^2$ (\revision{$\mathcal{E}_2^{rel}$}) and $L^\infty$ (\revision{$\mathcal{E}_\infty^{rel}$}) errors in $[0,15\,\textrm{m}]$ of the extended DG discretization with respect to a single-domain DG discretization, viscous Burgers' equation, sinusoidal initial data, constant scaling parameter $\beta=15$. $T=12\,\textrm{s}$, $\Delta t=10^{-2}\,\textrm{s}$.}
\centering
\begin{tabular}{cccc}\toprule\midrule
$N$ & $q$           & $\mathcal{E}_2^{\textrm{rel}}$ & $\mathcal{E}_\infty^{\textrm{rel}}$\\\midrule
\multirow{5}{*}{15} & 5 & 1.54E-02 & 3.88E-02 \\
                    & 10 & 1.42E-03 & 3.62E-03 \\
                    & 20  & 6.70E-04 & 1.71E-03 \\
                    & 40  & 6.63E-04 & 1.70E-03 \\
                    & 80  & 6.66E-04 & 1.71E-03 \\
\midrule
\multirow{5}{*}{30} & 5 & 8.27E-03 & 3.95E-02 \\
                    & 10 & 8.89E-04 & 4.29E-03 \\
                    & 30 & 1.38E-04 & 5.96E-04 \\
                    & 60 & 1.23E-04 & 5.21E-04 \\
                    & 100 & 1.20E-04 & 5.17E-04 \\
\bottomrule\end{tabular}\label{tab:burgers_wave_errors_const_beta}
\end{table}}
\revision{\begin{table}[htbp]\footnotesize\caption{Relative $L^2$ (\revision{$\mathcal{E}_2^{rel}$}) and $L^\infty$ (\revision{$\mathcal{E}_\infty^{rel}$}) errors in $[0,15\,\textrm{m}]$ of the extended DG discretization with respect to a single-domain DG discretization, viscous Burgers' equation, sinusoidal initial data, scaling parameter $\beta$ chosen for matching grid spacing at the finite/semi-infinite interface. $T=12\,\textrm{s}$, $\Delta t=10^{-2}\,\textrm{s}$.}
\centering
\begin{tabular}{ccccc}\toprule\midrule
$N$ & $q$ & $\beta$             & $\mathcal{E}_2^{\textrm{rel}}$ & $\mathcal{E}_\infty^{\textrm{rel}}$\\\midrule
\multirow{5}{*}{15} & 5 & 0.6 & 1.74E-01 & 4.96E-01 \\
                    & 10 & 0.35 & 5.03E-03 & 1.05E-02 \\
                    & 20 & 0.17 & 1.27E-02 & 3.54E-02 \\
                    & 40 & 0.09 & 9.92E-03 & 2.39E-02 \\
                    & 80 & 0.045 & 9.61E-03 & 2.33E-02 \\
\midrule
\multirow{5}{*}{30} & 5 & 1.2 & 8.74E-02 & 4.63E-01 \\
                    & 10 & 0.65 & 4.48E-02 & 2.02E-01 \\
                    & 30 & 0.24 & 1.71E-04 & 8.34E-04 \\
                    & 60 & 0.12 & 1.25E-03 & 6.63E-03 \\
                    & 100 & 0.075 & 2.36E-03 & 1.14E-02 \\
\bottomrule\end{tabular}\label{tab:burgers_wave_errors_matching_beta}
\end{table}}

\newpage

\subsection{Efficiency of the extended DG scheme in absorbing layer tests}
The second set of tests assesses the performance of the extended DG scheme in the absorption of perturbations leaving the finite subdomain when
an artificial damping term $-\gamma c, $ with  $\gamma\geq 0$, is added to the model equations' right-hand side on the the semi-infinite subdomain $[L,+\infty)$. As in \cite{benacchio:2013, benacchio:2019} we choose a sigmoid of the form 
\begin{equation}
    \gamma(z)=\frac{\Delta\gamma}{1+\text{exp}\left(\displaystyle\frac{\alpha L_0-z+L}{\sigma_D}\right)},
\end{equation}
where $\Delta\gamma$ is the sigmoid amplitude, $\alpha\in[0,1]$ the position of the sigmoid inside the absorbing layer, $L_0$ the spatial extension of the semi-infinite region, i.e. the distance between the first and the last Gauss-Laguerre-Radau nodes, and $\sigma_D$ the sigmoid steepness. As in \cite{benacchio:2013, benacchio:2019} we set $\alpha=0.3$ $m^{-1}$ and $\sigma_D=L_0/18$.

\subsubsection*{Advection-diffusion equation: Gaussian data}

In a first experiment, we consider the linear advection-diffusion equation \eqref{eq:advdifflin} and analyze the damping of a Gaussian profile defined by \eqref{eq:Gaussian} and initially placed inside the finite region $[0,L]$. To this end, we place the interface at $L=1000\,\textrm{m}$, set the initial data parameters $z_c=750\,\textrm{m}$, $\sigma_c=50\,\textrm{m}$, and  $\mu=1\,\textrm{m}^2/\textrm{s}$ and $u=1\,\textrm{m}/\textrm{s}$. The crest moves across the finite region, crosses the interface and is damped in the semi-infinite region. Spurious reflections into the finite region, measured as absolute errors of the computed solution in the finite region taking the absence of perturbation as reference, are below $10^{-3}$ for a range of values for the semi-infinite Laguerre modes $q$, finite subdomain subintervals $N$, and number of time steps $n$ -- and below $10^{-5}$ for the smallest $q=5$  (Table \ref{tab:dampGaussp2_fullyDG}). 
\begin{table}[htbp]\footnotesize
\caption{Absolute $L^2$ (\revision{$\mathcal{E}_2$}) and $L^\infty$ (\revision{$\mathcal{E}_\infty$}) residual errors in the finite region for the damping of a Gaussian perturbation with the extended DG scheme, linear advection-diffusion equation\revision{, several choices of number of modes $q$ and scaling parameter $\beta$ in the semi-infinite subdomain}. $C=0.33$, $p=2$, $T=500\,\textrm{s}$.}\label{tab:dampGaussp2_fullyDG}
\centering
\begin{tabular}{cccccc}\toprule\midrule
$q$                    & $N$ & $n$  & $\beta$ & $\mathcal{E}_2$ & $\mathcal{E}_\infty$ \\\midrule
40 & \multirow{5}{*}{400} & \multirow{5}{*}{600} & 1/28 & 9.22E-05 & 1.00E-04\\
30 &  & & 1/21 & 5.97E-06 & 6.75E-06\\
20 & & & 2/29 & 2.49E-05 & 2.61E-05\\
10 & & & 2/15 & 1.82E-06 & 1.24E-06\\
 5 & & & 1/4  & 1.51E-06 & 8.06E-07\\
\cmidrule(l{-.25pt}r{-.25pt}){1-6}
30 & \multirow{4}{*}{300} & \multirow{4}{*}{450} & 1/28 & 3.44E-04 & 3.42E-04\\
20 & & & 1/19 & 2.51E-04 & 2.51E-04\\
10 & & & 1/10 & 4.66E-06 & 4.25E-06\\
 5 & & & 11/60  & 1.65E-06 & 1.00E-06\\
\cmidrule(l{-.25pt}r{-.25pt}){1-6}
20 & \multirow{3}{*}{250} & \multirow{3}{*}{375} & 1/23 & 2.22E-04 & 2.05E-04\\
10 &  &  & 1/12 & 1.21E-05 & 1.10E-05\\
 5 & &  & 1/6  & 1.70E-06 & 1.07E-06\\
\cmidrule(l{-.25pt}r{-.25pt}){1-6}
10 & \multirow{2}{*}{200} & \multirow{2}{*}{300}& 1/15 & 4.25E-05 & 3.58E-05\\
 5 &  & & 1/7  & 1.86E-06 & 1.23E-06\\
\bottomrule\end{tabular}
\end{table}

\begin{figure}[h]
    \centering
    \includegraphics[width=0.5\textwidth]{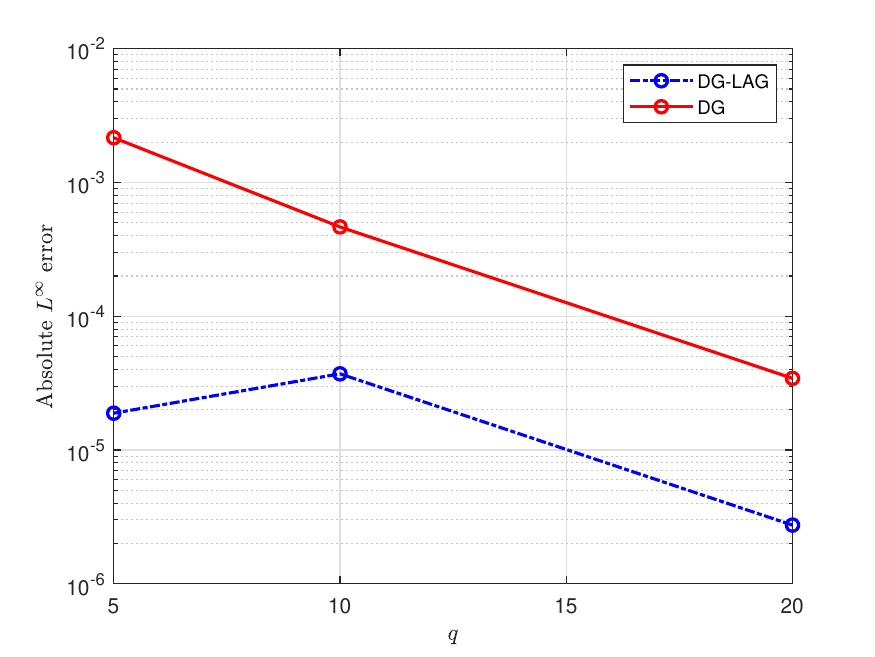}
    \caption{Damping of a Gaussian initial datum in the linear advection-diffusion equation. Absolute $L^\infty$ errors over the finite subdomain $[0,\,8\,\textrm{m}]$ of the extended DG solution (DG-LAG) and single-domain DG solution on a non-homogeneous grid (DG).}
    \label{fig:err_comp_abs_layer}
\end{figure}

\revision{\noindent In a second test with the same Gaussian initial data, we compare the extended DG scheme with a single-domain DG discretization in terms of efficiency of the absorbing layer implemented in the semi-infinite part subdomain. The interface in the coupled scheme is now placed at $L=8\,\textrm{m}$, the initial data parameters are $z_c=6\,\textrm{m}$ and $\sigma_c=1\,\textrm{m}$ and the physical parameters are $\mu=0.1\,\textrm{m}^2/\textrm{s}$ and $u=2\,\textrm{m}/\textrm{s}$. We run the simulation until the final time $T=4\,\textrm{s}$, with time step $\Delta t=0.02\,\textrm{s}$, and we choose $N=500$ intervals in $[0,L]$, with $q=20$, $10$ or $5$ Laguerre modes in $[L,+\infty)$. 

\noindent To make the absorbing layer comparison as fair as possible, we compute the single-domain DG solution on a non-uniform grid in $[0,L+L_0]$, such that the endpoints of the sub-intervals in $[L,L+L_0]$ coincide with the Laguerre nodes. By doing so, the single-domain DG grid is made of $N+q$ sub-intervals. Figure \ref{fig:abs_lay_comp} shows the two solutions at the final time for $q=20$. The small number of intervals in $[L+L_0]$ makes it challenging for the single-domain DG scheme to efficiently damp the outgoing signal. On the other hand, spectral accuracy in the same interval allows the Laguerre subdomain within the extended DG scheme to reduce the amplitude more significantly with the same number of modes. 

\noindent In terms of reflections into the finite domain, Table \ref{Tab:abs_lay_comp} shows the residual (absolute) errors in $[0,L]$ for the two schemes at the final time. Errors in the extended DG scheme are at least one order of magnitude smaller than the corresponding values in the single-domain DG case, and the difference becomes larger as $q$ decreases, reaching two orders of magnitude for $q=5$ (Figure \ref{fig:err_comp_abs_layer}). A small number of Laguerre modes is sufficient to obtain a good accuracy from the coupled scheme at significantly lower computational cost.

\begin{figure}[h]
    \centering
    \includegraphics[width=0.6\textwidth]{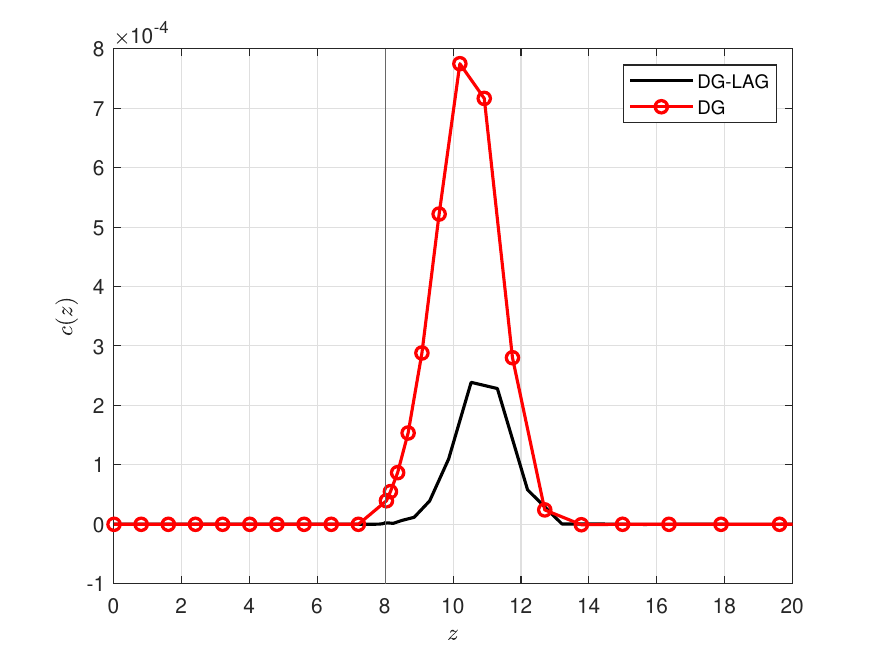}
    \caption{Damping of a Gaussian initial datum in the linear advection-diffusion equation. $T=4\,\textrm{s}$, $N=500$, $q=20$, $\Delta t=10^{-2}\,\textrm{s}$. Black line: extended DG scheme. Red line: Single-domain DG discretization on a non-uniform grid.}
    \label{fig:abs_lay_comp}
\end{figure}

\begin{table}[htbp]\footnotesize
\caption{Damping of a Gaussian initial datum in the linear advection-diffusion equation. Absolute $L^2$ (\revision{$\mathcal{E}_2$}) and $L^\infty$ (\revision{$\mathcal{E}_\infty$}) errors in $[0,L]$ of the extended DG solution (DG-LAG) and single-domain DG solution on a non-homogeneous grid (DG), for $q=20$, $q=10$ and $q=5$ modes in the semi-infinite subdomain. $L=8\,\textrm{m}$, $p=1$, $\mu=0.1\,\textrm{m}^2/\textrm{s}$, $u=2\,\textrm{m}/\textrm{s}$, $T=4\,\textrm{s}$, $\Delta t=0.02\,\textrm{s}$.}\label{Tab:abs_lay_comp}
\centering
\begin{tabular}{cccc}\toprule\midrule
$q$                 &        & $\mathcal{E}_2$ & $\mathcal{E}_\infty$ \\ \cmidrule(l{-.25pt}r{-.25pt}){1-4}
\multirow{2}{*}{20} & DG-LAG & 5.56E-07 & 2.74E-06 \\
                    & DG     & 6.09E-06 & 3.43E-05 \\ \cmidrule(l{-.25pt}r{-.25pt}){1-4}
\multirow{2}{*}{10} & DG-LAG & 5.80E-06 & 3.71E-05 \\
                    & DG     & 8.03E-05 & 4.65E-04 \\ \cmidrule(l{-.25pt}r{-.25pt}){1-4}
\multirow{2}{*}{5}  & DG-LAG & 2.31E-06 & 1.88E-05 \\
                    & DG     & 3.69E-04 & 2.16E-03 \\
\bottomrule\end{tabular}
\end{table}
\noindent
}

\newpage

\subsubsection*{Advection-diffusion equation: Wave train}

\begin{figure}
\centering
\includegraphics[width=.45\textwidth]{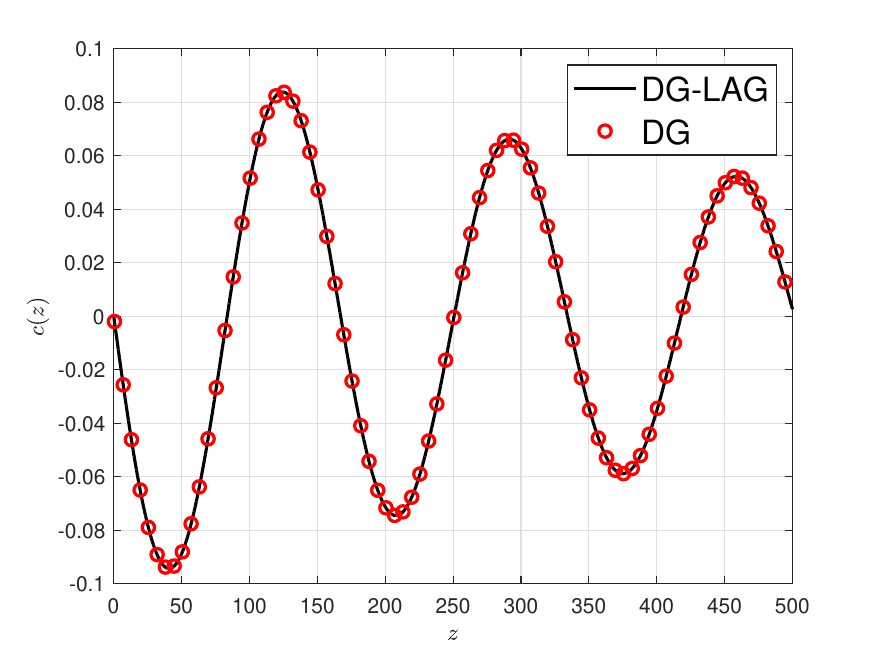}\includegraphics[width=.45\textwidth]{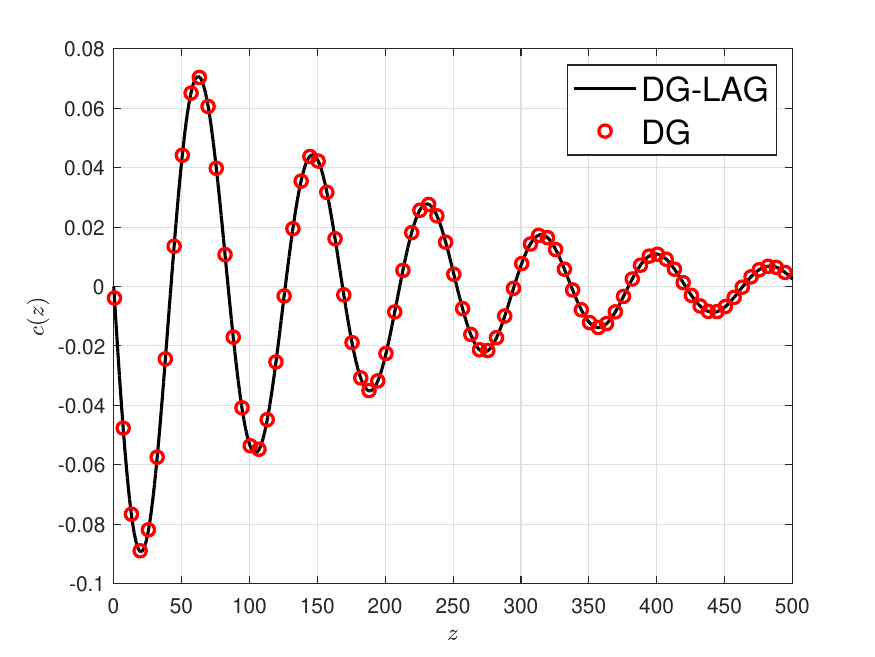}
\caption{Damping of a wave train with wavenumber $k=30$ (left), $k=60$ (right), linear advection-diffusion equation. Solid line: extended DG scheme. Red circles: single-domain DG scheme. $A=0.1\,\textrm{m}$, $q=30$, $N=600$, $\beta=0.143$, $T=5000\,\textrm{s}$, $n=16000$.}
    \label{fig:dampwave1}
\end{figure}

\noindent Next, we consider a wave train case, obtained by imposing a Dirichlet boundary condition 
\begin{equation}
c(0,t)=A\text{sin}(2\pi k/T t)
\end{equation}
at the left endpoint $z=0$. The initial condition is $c_0=0$. The wave train is generated at $z=0$, crosses the finite region $[0,L]$ and is damped by the absorbing layer, where we set $\Delta\gamma=2A$. Numerical parameters are set as $L=500\,\textrm{m}$, $\mu=1\,\textrm{m}^2/\textrm{s}$, $u=1\,\textrm{m}/\textrm{s}$, $T=5000\,\textrm{s}$, and $n=16000$ time steps.
On a range of choices for the wave number, amplitude, and Laguerre modes $q$, the extended DG scheme absorbs outgoing perturbations with relative errors computed in $[0,L]$ of less than $10^{-4}$ for $q=15$, and at most $1.3\times10^{-3}$ for $q=5$, with respect to a reference single-domain DG solution on $[0,2L]$ (Figure \ref{fig:dampwave1} and Table \ref{Tab:dampwave}). Results are comparable with those obtained in \cite{benacchio:2019} for the inviscid shallow water system with a different coupling approach. The efficiency of the tool is competitive -- for $q=5$, less than a hundredth of the computational cost for $N=600$, and less than five thousandths of the computational cost for $N=1200$, is due to the absorbing layer.  

\begin{table}[htbp]\footnotesize
\caption{Damping of a wave train, linear advection-diffusion equation. Relative $L^2$ (\revision{$\mathcal{E}^{rel}_2$}) and $L^\infty$ (\revision{$\mathcal{E}^{rel}_\infty$}) errors in $[0,L]$ of the extended DG solution with respect to a single-domain DG solution, \revision{for $q=15$ and $q=5$} modes in the semi-infinite subdomain\revision{, and several choices of amplitude $A$, wavenumber $k$, and elements $N$ in the finite subdomain.} $L=500\,\textrm{m}$, $p=1$, $\mu=1\,\textrm{m}^2/\textrm{s}$, $u=1\,\textrm{m}/\textrm{s}$, $T=5000\,\textrm{s}$, $n=16000$.}\label{Tab:dampwave}
\centering
\begin{tabular}{ccccccc}\toprule\midrule
$\mathbf{q}$ & $A$                    & $k$ & $N$  & $\beta$ & $\mathcal{E}_2^{\textrm{rel}}$ & $\mathcal{E}_\infty^{\textrm{rel}}$\\\cmidrule(l{-.25pt}r{-.25pt}){1-7}
& \multirow{2}{*}{0.025} & 30  & 600    & 0.286  & 1.60E-06 & 2.14E-05\\
& & 60  & 1200   & 0.571  & 1.66E-07 & 2.24E-06\\
\cmidrule(l{-.25pt}r{-.25pt}){2-7}
\multirow{2}{*}{$\mathbf {15}$}&\multirow{2}{*}{0.05}  & 30  & 600    & 0.286  & 2.26E-06 & 2.99E-05\\
& & 60  & 1200   & 0.571  & 2.61E-07 & 3.30E-06\\
\cmidrule(l{-.25pt}r{-.25pt}){2-7}
&\multirow{2}{*}{0.1}  & 30  & 600    & 0.286  & 2.49E-06 & 3.13E-05\\
& & 60  & 1200   & 0.571  & 4.76E-07 & 6.02E-06\\\midrule
 & \multirow{2}{*}{0.025} & 30  & 600    & 0.74  & 7.56E-05 & 1.04E-03\\
& & 60  & 1200   & 1.48  & 4.27E-06 & 3.02E-05\\
\cmidrule(l{-.25pt}r{-.25pt}){2-7}
\multirow{2}{*}{$\mathbf {5}$} & \multirow{2}{*}{0.05}  & 30  & 600    & 0.74  & 3.70E-05 & 5.13E-04\\
&  & 60  & 1200   & 1.48  & 7.34E-06 & 5.49E-05\\
\cmidrule(l{-.25pt}r{-.25pt}){2-7}
&\multirow{2}{*}{0.01}  & 30  & 600    & 0.74  & 3.10E-05 & 4.32E-04\\
&  & 60  & 1200   & 1.48  & 1.14E-05 & 8.58E-05\\
\bottomrule
\bottomrule\end{tabular}
\end{table}

\newpage

\subsubsection*{Burgers' equation: Gaussian initial data}

Finally, we consider the Burgers' equation. We place the interface at $L=30\,\textrm{m}$, and center an initial Gaussian profile inside the bounded region ($z_c=25\,\textrm{m}$, $\sigma_c=1\,\textrm{m}$). We run the extended DG scheme until $T=3600\,\textrm{s}$, when most of the initial perturbation has left the bounded region. Residual errors in the finite region $[0,L]$ with the extended DG scheme with respect to a single-domain DG solution on $[0,10L/3]$ are below one percent for as few as $q=5$ modes in the semi-infinite subdomain (Table \ref{Tab:burgers_inside}). Qualitatively equivalent results are obtained when placing the initial data at the interface, $z_c=30\,\textrm{m}$ (not shown).

\begin{table}[htbp]\footnotesize\caption{Damping of a Gaussian profile initially centred at $z_c=25\,\textrm{m}$, viscous Burgers' equation. Relative $L^2$ (\revision{$\mathcal{E}^{rel}_2$}) and $L^\infty$ (\revision{$\mathcal{E}^{rel}_\infty$}) errors in $[0,L]$ of the extended DG scheme compared with a single-domain DG solution, several choices of the number of modes $q$ in the semi-infinite subdomain.  $T=3600\,\textrm{s}$, $\Delta t = 0.1\,\textrm{s}$, $N=30$, $p=1$, $\Delta\gamma=2$, $\sigma_c=1\,\textrm{m}$.}\label{Tab:burgers_inside}
\centering
\begin{tabular}{cccc}\toprule\midrule
$q$ & $\beta$             & $\mathcal{E}^{rel}_2$ & $\mathcal{E}^{rel}_\infty$\\\midrule
60  & 0.06 & 2.12E-03 & 2.06E-03 \\
40  & 0.09 & 2.13E-03 & 2.06E-03 \\
20  & 0.175 & 2.15E-03 & 2.09E-03 \\
10  & 0.34 & 2.39E-03 & 2.31E-03 \\
5   & 0.68 & 7.13E-03 & 6.82E-03 \\
\bottomrule\end{tabular}
\end{table}

 \section{Conclusions and perspectives}
\label{sec:conclu} \indent  

This paper proposed an extended DG approach for the numerical simulation of nonlinear \revision{advection-diffusion} problems on unbounded domains. Built on earlier developments of coupled DG-Laguerre discretizations for purely hyperbolic systems, the scheme models a finite portion of the semi-infinite half-line using standard Legendre basis functions and the adjacent unbounded portion using scaled Laguerre basis functions. 

Compared to a standard DG discretization, the extended DG scheme only differs for the presence of two off-diagonal terms in the system matrix, representing the numerical fluxes at the finite/semi-infinite interface. The new framework improves on previous endeavours that used bespoke coupling strategies, and provides a completely seamless coupling approach.

\revision{The stability of the extended DG scheme was analyzed in the linear case, proving that the resulting matrix has eigenvalues with negative real part in the inviscid case. In the viscous case, the same analysis was performed empirically and
 the resulting matrix was found to have eigenvalues with negative real part in typical configurations and independently of the P\'eclet number.} The analysis and numerical experiments used scaled Laguerre basis functions, Gauss-Laguerre-Radau quadrature in the unbounded subdomain, and Dirichlet boundary conditions.  Results covering other possible options are reported in \cite{vismara:2020} and are summarized in the Appendix, corroborating the findings for the standalone Laguerre scheme in a purely hyperbolic framework \cite{benacchio:2019}. 
To the best of the authors' knowledge, a stability analysis on numerical schemes using different sets of basis functions as the one presented in this work is not currently available in the literature.

The correctness of the extended DG scheme, particularly regarding the finite/semi-infinite interface fluxes, was validated in a series of numerical experiments with the linear homogeneous and non-homogeneous advection-diffusion equation and the nonlinear viscous Burgers' equation. By comparison with a standard single domain implementation, spurious signals due to the presence of different basis functions are of negligible entity on a range of spatial resolutions, thereby complementing and strengthening results obtained with hyperbolic systems in \cite{benacchio:2013, benacchio:2019}.

In tests where the semi-infinite subdomain featured a reactive damping term, the extended DG scheme displayed compelling performance in efficiently absorbing outgoing waves in linear and nonlinear models. A very small number of Laguerre modes, both in absolute terms and as a proportion of the total computational load, was sufficient to damp single Gaussian signals and wave trains without spurious phenomena spoiling the simulation in the finite subdomain.  \revision{Notably, the extended DG scheme displayed reflections in the finite subdomain with maximum amplitude more than one order of magnitude smaller compared with a single-domain DG scheme using a non-uniform grid and the same number of modes. In addition, the advantage using the proposed extended scheme grows with decreasing number of modes in the semi-infinite subdomain.  While the results were obtained with a linear advection-diffusion model, we expect these findings to be corroborated on nonlinear systems and tests with more complex wave dynamics, making the extended DG scheme an interesting technique for the discretization of fluid dynamics problems on unbounded domains.}

The results achieved in this work offer a number of perspectives for future investigation. First, a similar extended DG approach can be developed coupling a strong form, nodal DG discretization on the finite domain to the strong form, nodal approach with scaled Laguerre functions for the semi-infinite domain, a choice that displayed stability advantages in the large  P\'eclet number case. The scheme can then be implemented in multiple dimensions, using tensor product-based discretization approaches on semi-infinite strips or circular domains, where the problem is discretized using the extended DG scheme in the vertical or radial direction and a discontinuous Galerkin approach in the horizontal or azimuthal direction. Such a model may find applications, for example, in the modelling of the solar corona. The extension to systems of parabolic equations or to non-linear diffusion may be considered, such as are found in turbulence modelling. From a more theoretical perspective, the possibility to prove inf-sup conditions for the extended DG approach could also be investigated.  

\section*{Acknowledgements}
This work summarizes and extends results obtained by the first author (F.V.) in his Master's Thesis in Mathematical Engineering \cite{vismara:2020}, discussed at Politecnico di Milano in 2020 and prepared under the supervision of the other two authors.
T.B. and L.B. have been supported by the ESCAPE-2 project of the Horizon 2020 research and innovation programme (grant agreement No 800897). Two anonymous reviewers are gratefully acknowledged for their critical comments, which have helped to improve the presentation of the paper's results.

\section*{Data availability}

The datasets generated during the current study are available from the corresponding author on reasonable request.

\section*{Conflict of interest}
The authors declare that they have no conflict of interest.

\newpage

\appendix
\section{Alternative discretizations on the semi-infinite subdomain}
\label{sec:alternatives} \indent  
We summarize here the results presented in \cite{vismara:2020} on the analysis of various Laguerre-based discretizations of the advection-diffusion equation with constant coefficients on ${\mathbb R}^+=[0,+\infty)$. 
For the purpose of deriving some discretizations, it can be helpful to reformulate equation \eqref{eq:advdifflin}, which we report here for convenience,
\begin{equation}
    \frac{\partial c}{\partial t}+u\frac{\partial c}{\partial z}=\mu\frac{\partial^2 c}{\partial z^2} 
\end{equation}
as a system of first order equations
\begin{align}
\label{eq:advdiff1st}
\begin{split}
 &   \frac{\partial c}{\partial t} - \mu \frac{\partial v}{\partial z} + uv =0  \\ 
 &   \frac{\partial c}{\partial z} - v = 0. 
 \end{split}
\end{align}
 We assume that solutions vanish at infinity
\begin{equation}
\label{eq:vanish}
    \lim_{z\rightarrow+\infty}c(z,t)=0
\end{equation}
and that either Dirichlet boundary conditions
  \begin{equation}
  \label{eq:dirichlet}
 c(0,t) = c_L
\end{equation}
  or Neumann boundary conditions
  \begin{equation}
  \label{eq:neumann}
    \frac{\partial c}{\partial z}(0,t) = Dc_L
\end{equation}
are applied at $z=0.$  We require that $\mu>0$ (ellipticity condition) and $u>0$. In this case, the Dirichlet datum at $z=0$ corresponds to an inflow boundary condition, which guarantees well-posedness for the hyperbolic part. We analyze several possible space discretizations, in order to determine which one shows the best stability properties and can therefore be chosen for the extended DG scheme in conjunction with the Legendre basis in the finite sub-domain. As done in \cite{benacchio:2019} for the pure advection problem, we discretize the PDE system \eqref{eq:advdiff1st} in space, obtaining, after substitution of the discretization of the second equation in \eqref{eq:advdiff1st} into the first, a system of ordinary differential equations of the form
\begin{equation}
\frac{d\mathbf{c}}{dt}=\mathbf{A}\mathbf{c}+\mathbf{g},
\end{equation}
where $\mathbf{c}$ is the unknown vector of the expansion of the solution and $\mathbf{g}$ contains the contribution of boundary conditions at $z=0$, and we study the eigenvalue structure of the matrix $\mathbf{A}$. The corresponding discretization scheme is stable if all the eigenvalues have non-positive real part. \\ We analyse the following discretizations:
\begin{itemize}

\item \emph{Weak form}. We multiply  \eqref{eq:advdiff1st} by a test function, integrate by parts and use either Gauss-Laguerre-Radau (GLR) or Gauss-Laguerre (GL) quadrature rules.
Two different approaches are possible. In a modal approach, entries of the unknown vector $\mathbf{c}$ are the coefficients of the expansion of the solution in the orthogonal basis of Laguerre functions or Laguerre polynomials. In a nodal approach, the basis functions are the Lagrange basis functions associated with the integration nodes, so that the unknown vector contains the nodal values of the approximate solution. Furthermore, the numerical solution can be expanded in a basis of either scaled Laguerre functions or scaled Laguerre polynomials.

\item \emph{Strong form}. In this case we directly discretize the strong formulation \eqref{eq:advdiff1st} using a collocation approach and GLR quadrature rules. This is the only practical choice if Dirichlet boundary conditions have to be imposed, because the GLR nodes include the left endpoint of the semi-infinite subdomain, unlike the GL nodes.
\end{itemize}

\noindent We now summarize some definitions we need to introduce the different variants of the matrix $\mathbf{A}$ and vector $\mathbf{g}.$ For discretizations based on Laguerre functions, we define the matrix $\mathbf{\hat{L}}=\{\hat{l}_{ij}\}$ 
with entries  such that
\begin{equation}\label{eq:matrixLf}
    \hat{l}_{ij}=\begin{cases}
    1/2 & i=j \\
    1 & j<i \\
    0 & j>i.
    \end{cases}
\end{equation}
If discretizations based on Laguerre polynomials are considered, we use the matrix $\mathbf{L}=\{l_{ij}\}$ defined~ as

\begin{equation}\label{eq:matrixLp}
    l_{ij}=\begin{cases}
    0 & i=j \\
    1 & j<i \\
    0 & j>i
    \end{cases}
\end{equation}
For nodal discretizations based on the weak form and on scaled Laguerre functions, we then denote by $z_j^\beta$ the $j$-th GLR or GL quadrature node, by $h_j^\beta(z)$ the associated Lagrangian polynomial, by $\omega_i$ the $i$-th quadrature weight, and by $\hat{d}_{ij}^\beta$ the entries of the GLR or GL differentiation matrix $\mathbf{\hat{D}}_\beta$ associated with scaled Laguerre functions, defined as follows:
\begin{itemize}
    \item GL nodes
    \begin{equation} \hat{d}^\beta_{ij} = \begin{cases} 
          \dfrac{\hat{\mathscr L}_q^\beta(z_i^\beta)}{(z_i^\beta-z_j^\beta)\hat{\mathscr L}_q^\beta(z_j^\beta)} & i\neq j \\ \\
          -\dfrac{q+2}{2z_i^\beta} & i=j \\
       \end{cases}
    \end{equation}
    \item GLR nodes
    \begin{equation} \hat{d}^\beta_{ij} = \begin{cases} 
          \dfrac{\hat{\mathscr L}_{q+1}^\beta(z_i^\beta)}{(z_i^\beta-z_j^\beta)\hat{\mathscr L}_{q+1}^\beta(z_j^\beta)} & i\neq j \\ \\
          0 & i=j\neq 0 \\ \\
          -\beta\dfrac{q+1}{2} & i=j=0. \\
       \end{cases}
    \end{equation}
\end{itemize}
We also define as $\boldsymbol{\hat{\Omega}}_\beta$ the diagonal matrix with the quadrature weights $\hat{\omega}_i^\beta$ on the diagonal.
For a nodal discretization based on Laguerre polynomials, instead, the differentiation matrix $\mathbf{D}_\beta$  has entries $d^\beta_{ij}$  defined as:
\begin{itemize}
    \item GL nodes
    \begin{equation} {d}_{ij}^\beta = \begin{cases} 
          \dfrac{{\mathscr L}_q^\beta(z_i^\beta)}{(z_i^\beta-z_j^\beta){\mathscr L}_q^\beta(z_j^\beta)} & i\neq j \\ \\
          \dfrac{\beta z_i^\beta-q-2}{2z_i^\beta} & i=j \\
       \end{cases}
    \end{equation}
    \item GLR nodes
    \begin{equation} {d}_{ij}^\beta = \begin{cases} 
          \dfrac{{\mathscr L}_{q+1}^\beta(z_i^\beta)}{(z_i^\beta-z_j^\beta){\mathscr L}_{q+1}^\beta(z_j^\beta)} & i\neq j \\ \\
          \dfrac{\beta}{2} & i=j\neq 0 \\ \\
          -\beta\dfrac{q}{2} & i=j=0 \\
       \end{cases}
    \end{equation}
\end{itemize}
We also set
\begin{gather*}
\mathbf{\hat{g}_1}=[(\hat{h}_0^\beta)^{\prime\prime}(z_1),\dots,(\hat{h}_0^\beta)^{\prime\prime}(z_q)],\;\mathbf{\hat{g}_2}=[(\hat{h}_0^\beta)^{\prime}(z_1),\dots,(\hat{h}_0^\beta)^{\prime}(z_q)],\\
\mathbf{g_1}=[(h_0^\beta)^{\prime\prime}(z_1),\dots,(h_0^\beta)^{\prime\prime}(z_q)],\;\mathbf{g_2}=[(h_0^\beta)^{\prime}(z_1),\dots,(h_0^\beta)^{\prime}(z_q)],\\
\mathbf{\hat{h}}=[\hat{h}^\beta_0(0),\dots,\hat{h}^\beta_q(0)],\; \mathbf{h}=[h^\beta_0(0),\dots,h^\beta_q(0)],\\
\mathbf{\widehat{W}}=\boldsymbol{\hat{\Omega}}_\beta^{-1}\mathbf{\hat{D}}_\beta^T\boldsymbol{\hat{\Omega}}_\beta,\; \mathbf{W}=\boldsymbol{\Omega}_\beta^{-1}\mathbf{D}_\beta^T\boldsymbol{\Omega}_\beta,\\
\mathbf{\hat{r}}=\boldsymbol{\hat{\Omega}}_\beta^{-1}\mathbf{\hat{h}},\; \mathbf{r}=\boldsymbol{\Omega}_\beta^{-1}\mathbf{h},\;\mathbf{e}=[1,\dots,1]^T\in\mathbf{R}^{q+1}.
\end{gather*} 
We also denote by $(\mathbf{\hat{D}}_\beta)_q$ for scaled Laguerre functions, and by $(\mathbf{D}_\beta)_q$ for scaled Laguerre polynomials, the matrices obtained from the differentiation matrices $\mathbf{\hat{D}}_\beta$ and $\mathbf{D}_\beta$ by removing the first row and the first column. Finally we denote by $(\mathbf{\hat{D}}_\beta)_0$ for scaled Laguerre functions, and $(\mathbf{D}_\beta)_0$ for scaled Laguerre polynomials, the matrices obtained from $\mathbf{\hat{D}}_\beta$ and $\mathbf{D}_\beta$ by replacing the first row with zeros. The expressions of matrix $\mathbf{A}$ and right-hand side $\mathbf{g}$ for the derived discretizations are summarized in Table \ref{tab:matrixA} -- note the two use of the matrix $\mathbf{\hat{L}}$  \eqref{eq:matrixLf} for scaled Laguerre functions and $\mathbf{L}$ \eqref{eq:matrixLp} for scaled Laguerre polynomials.

\begin{table}[htbp]\footnotesize
\centering
\caption{Definition of matrix $\mathbf{A}$ and vector $\mathbf{g}$ in the Laguerre discretization of \eqref{eq:advdifflin} or \eqref{eq:advdiff1st} for several formulations `Form', basis functions `BF' and boundary conditions `BC'. `Coll': collocation, `Nod': nodal, `Mod': modal, `Dir': Dirichlet, `Neu': Neumann, `LF': Scaled Laguerre Functions, `LP': Scaled Laguerre Polynomials. See text for symbol definitions.}
\begin{tabular}{lllll}
\toprule\midrule
Form & BF & BC & $\mathbf{A}$ & $\mathbf{g}$\\
\midrule
Coll & LF & Dir & $\mu(\mathbf{\hat{D}}_\beta^2)_q-u(\mathbf{\hat{D}}_\beta)_q$ & $\mu c_L\mathbf{\hat{g}_1}-uc_L\mathbf{\hat{g}_2}$  \\
Coll & LF & Neu & $\mu\mathbf{\hat{D}}_\beta(\mathbf{\hat{D}}_\beta)_0-u(\mathbf{\hat{D}}_\beta)_0$ & $\mu Dc_L\mathbf{\hat{g}_2}-uDc_L\mathbf{e_1}$\\
Coll & LP & Dir & $\mu(\mathbf{{D}}_\beta^2)_q-u(\mathbf{{D}}_\beta)_q$ & $\mu c_L\mathbf{g_1}-uc_L\mathbf{g_2}$ \\
Coll & LP & Neu & $\mu\mathbf{{D}}_\beta(\mathbf{{D}}_\beta)_0-u(\mathbf{{D}}_\beta)_0$ & $\mu Dc_L\mathbf{g_2}-uDc_L\mathbf{e_1}$\\ \midrule
Nod & LF & Dir  & $-\mu\mathbf{\hat{D}}_\beta\mathbf{\hat{W}} + u\mathbf{\hat{W}}$ & $-\mu c_L\mathbf{\hat{D}}_\beta\mathbf{\hat{r}}+ uc_L\mathbf{\hat{r}} $ \\
Nod & LF & Neu  & $-\mu\boldsymbol{\hat{\Omega}_\beta}^{-1}\mathbf{\hat{D}}_\beta^T\boldsymbol{\hat{\Omega}}_\beta\mathbf{\hat{D}}_\beta-u\mathbf{\hat{D}}_\beta$  & $-\mu Dc_L\mathbf{\hat{r}}$\\ 
Nod & LP & Dir  & $-\mu\mathbf{D}_\beta\mathbf{W}+\mu\beta\mathbf{D}_\beta+u\mathbf{W}-u\beta\mathbf{I}$ & $-\mu c_L\mathbf{D}_\beta\mathbf{r}+uc_L\mathbf{r}$\\ 
Nod & LP & Neu  & $-\mu\mathbf{W}\mathbf{D}_\beta+\mu\beta\mathbf{D}_\beta-u\mathbf{D}_\beta$  & $-\mu Dc_L\mathbf{r}$ \\ \midrule
Mod & LF & Dir  & $-\mu\beta^2\mathbf{\hat{L}}^T\mathbf{\hat{L}} - u\beta\mathbf{\hat{L}}$  & $\mu \beta^2c_L\mathbf{\hat{L}}^T\mathbf{e} + u\beta c_L\mathbf{e}$\\
Mod & LF & Neu  & $-\mu\beta^2\mathbf{\hat{L}}\mathbf{\hat{L}}^T + u\beta \mathbf{\hat{L}}^T$ & $-\mu \beta Dc_L \mathbf{e}$\\ 
Mod & LP & Dir  & $-\mu\beta^2\mathbf{L}^T(\mathbf{L}+\mathbf{I})-u\beta(\mathbf{L}+\mathbf{I})$ & $u\beta c_L\mathbf{e}$\\
Mod & LP & Neu & $-\mu\beta^2(\mathbf{L}+\mathbf{I})\mathbf{L}^T + u\beta \mathbf{L}^T$ & $-\mu \beta Dc_L \mathbf{e}$\\
\bottomrule
\end{tabular}\label{tab:matrixA}
\end{table}

As customary for the advection-diffusion problem, the stability property can be a function of the P\'eclet number, which is usually defined as $Pe=u\mathcal{L}/\mu$, where $\mathcal{L}$ is a reference length scale. For simplicity we choose the length scale $\mathcal{L}=1$, set $\mu=1$ and analyze the stability of $\mathbf{A}$ for a fixed value of $Pe$; the corresponding ranges for $\beta$ are shown in Table \ref{tab:pe1} for both scaled Laguerre functions and polynomials.

\begin{table}[htbp]\footnotesize\caption{Stability of $\mathbf{A}$ as a function of $\beta$: condition under which the largest real part of the eigenvalues is non-positive. $q=50$, $\mu=1$. `Neu': Neumann b.c., `Dir': Dirichlet b.c., `LF': Scaled Laguerre Functions, `LP': Scaled Laguerre Polynomials.}
\centering
\begin{tabular}{lllcccc}\toprule\midrule
& & & \multicolumn{2}{c}{\qquad LF} &\multicolumn{2}{c}{\qquad LP} \\
\cmidrule(l{5pt}r{-.5pt}){4-7}
&  & & Neu & Dir & Neu & Dir \\
\cmidrule(l{5pt}r{-.5pt}){4-7}
Strong & & & $\forall\beta$ & $\forall\beta$ & $\beta\leq 2.6Pe$ & $\beta\leq 3Pe$ \\
\cmidrule(l{5pt}r{-.5pt}){1-7}
\multirow{3}{*}{Weak} & \multirow{2}{*}{Nodal} & GLR  & $\beta\geq 0.58Pe$ & $\forall\beta$ & $0.017Pe\leq \beta\leq 2.83Pe$      & $\beta\leq 3Pe$ \\
                           &                        & GL & $\beta\geq2Pe$ & $\forall\beta$ & $0.25Pe\leq\beta\leq 2Pe$ & $\beta\leq 8.5Pe$ \\     
\cmidrule(lr{-.5pt}){2-7}
                           &                 Modal  &      & $\beta\geq0.58Pe$ & $\forall\beta$ & $0.017Pe\leq \beta\leq 2.83Pe$     & $\beta\leq 3Pe$                    \\
\bottomrule\end{tabular}\label{tab:pe1}
\end{table}
It can be observed that only the strong form discretizations based on Laguerre functions are stable for all boundary conditions and independently of the value of the P\'eclet number. Other discretizations based on Laguerre functions are instead stable under mild conditions on the value of $\beta$ as a function of the P\'eclet number. These conditions  become problematic only in the very large P\'eclet number limit. 

In this paper, only the weak form modal discretization based on Laguerre functions was considered for the extended DG scheme, due to its hierarchical nature, that allows in principle for an easy (and if necessary, dynamic) adjustment of the number of basis functions to perform  $p-$adaptation. The strong form nodal discretization based on Laguerre functions seems otherwise the most robust option and will be further studied as a basis for extended DG approaches in future work.  Discretizations based on Laguerre polynomials are instead only stable under more restrictive conditions, which also affect the choice of  $\beta$ in the small  P\'eclet number case.
These conclusions complement the results in \cite{benacchio:2019}, where the pure advection problem was discussed.  Such an analysis does not seem to have been carried out in the literature, to the best of the authors' knowledge.

%

\bibliographystyle{spmpsci}      
\bibliography{dg_laguerre_1d_parabolic_jsc}   

\begin{thebibliography}{10}
\providecommand{\url}[1]{{#1}}
\providecommand{\urlprefix}{URL }
\expandafter\ifx\csname urlstyle\endcsname\relax
  \providecommand{\doi}[1]{DOI~\discretionary{}{}{}#1}\else
  \providecommand{\doi}{DOI~\discretionary{}{}{}\begingroup
  \urlstyle{rm}\Url}\fi

\bibitem{akmaev:2011}
Akmaev, R.: Whole atmosphere modeling: Connecting terrestrial and space
  weather.
\newblock Reviews of Geophysics \textbf{49} (2011)

\bibitem{appelo:2009}
Appel{\"o}, D., Colonius, T.: A high-order super-grid-scale absorbing layer and
  its application to linear hyperbolic systems.
\newblock Journal of Computational Physics \textbf{228}(11), 4200--4217 (2009)

\bibitem{arnold:1982}
Arnold, D.: An interior penalty finite element method with discontinuous
  elements.
\newblock SIAM Journal of Numerical Analysis \textbf{19}, 742--760 (1982)

\bibitem{arnold:2002}
Arnold, D., Brezzi, F., Cockburn, B., Marini, L.: Unified analysis of
  discontinuous {G}alerkin methods for elliptic problems.
\newblock SIAM Journal of Numerical Analysis \textbf{39}, 1749--1779 (2002)

\bibitem{astley:2000}
Astley, R.: Infinite elements for wave problems: a review of current
  formulations and an assessment of accuracy.
\newblock International Journal of Numerical Methods in Engineering
  \textbf{49}(7), 951--976 (2000)

\bibitem{benacchio:2013}
Benacchio, T., Bonaventura, L.: Absorbing boundary conditions: a spectral
  collocation approach.
\newblock International Journal of Numerical Methods in Fluids \textbf{72}(9),
  913--936 (2013).
\newblock \doi{10.1002/fld.3768}.
\newblock \urlprefix\url{http://dx.doi.org/10.1002/fld.3768}

\bibitem{benacchio:2019}
Benacchio, T., Bonaventura, L.: An extension of {DG} methods for hyperbolic
  problems to one-dimensional semi-infinite domains.
\newblock Applied Mathematics and Computation \textbf{350}, 266--282 (2019)

\bibitem{black:1998}
Black, K.: Spectral elements on infinite domains.
\newblock SIAM Journal of Scientific Computing \textbf{19}, 1667--1681 (1998)

\bibitem{bonaventura:2000}
Bonaventura, L.: A {S}emi--implicit {S}emi--{L}agrangian {S}cheme {U}sing the
  {H}eight {C}oordinate for a {N}onhydrostatic and {F}ully {E}lastic {M}odel of
  {A}tmospheric {F}lows.
\newblock Journal of Computational Physics \textbf{158}(2), 186--213 (2000)

\bibitem{bonaventura:2017}
Bonaventura, L., Fern\'andez-Nieto, E., Garres-D\'iaz, J., Narbona-Reina, G.:
  Multilayer shallow water models with locally variable number of layers and
  semi-implicit time discretization.
\newblock Journal of Computational Physics \textbf{364}, 209--234 (2017)

\bibitem{dea:2011}
Dea, J.: An experimental adaptation of {H}igdon-type non-reflecting boundary
  conditions to linear first-order systems.
\newblock Journal of Computational and Applied Mathematics \textbf{235},
  1354--1366 (2011)

\bibitem{engquist:1977}
Engquist, B., Majda, A.: Absorbing boundary conditions for numerical simulation
  of waves.
\newblock Mathematics of Computation \textbf{31}(139), 629--651 (1977)

\bibitem{gerdes:2000}
Gerdes, K.: A review of infinite element methods for exterior {H}elmholtz
  problems.
\newblock Journal of Computational Acoustics \textbf{8}(1), 43--62 (2000)

\bibitem{giraldo:2013}
Giraldo, F., Kelly, J., Constantinescu, E.: Implicit-explicit formulations of a
  three-dimensional nonhydrostatic unified model of the atmosphere ({NUMA}).
\newblock SIAM Journal on Scientific Computing \textbf{35} (2013)

\bibitem{gu:2021}
Gu, D., Wang, Z.: Orthogonal jacobi rational functions and spectral methods on
  the half line.
\newblock Journal of Scientific Computing \textbf{88}(1), 1--27 (2021)

\bibitem{israeli:1981}
Israeli, M., Orszag, S.: Approximation of {R}adiation {B}oundary {C}onditions.
\newblock Journal of Computational Physics \textbf{41}, 115--135 (1981)

\bibitem{jackson:2019}
Jackson, D., Fuller-Rowell, T., Griffin, D., Griffith, M., Kelly, C., Marsh,
  D., Walach, M.: Future directions for whole atmosphere modeling: Developments
  in the context of space weather.
\newblock Space Weather \textbf{17}, 1342--1350 (2019)

\bibitem{klemp:1983}
Klemp, J., Durran, D.: An {U}pper {B}oundary {C}ondition {P}ermitting
  {I}nternal {G}ravity {W}ave {R}adiation in {N}umerical {M}esoscale {M}odels.
\newblock Journal of Atmospheric Sciences \textbf{111}, 430--444 (1983)

\bibitem{klemp:1978}
Klemp, J., Lilly, D.: Numerical {S}imulation of {H}ydrostatic {M}ountain
  {W}aves.
\newblock Journal of Atmospheric Sciences \textbf{35}, 78--107 (1978)

\bibitem{rasch:1986}
Rasch, P.: Toward atmospheres without tops: {A}bsorbing upper boundary
  conditions for numerical models.
\newblock Quarterly Journal of the Royal Meteorological Society \textbf{112},
  1195--1218 (1986)

\bibitem{riviere:2008}
Rivière, B.: Discontinuous Galerkin Methods for Solving Elliptic and Parabolic
  Equations: Theory and Implementation.
\newblock SIAM (2008)

\bibitem{shen:2001}
Shen, J.: Stable and efficient spectral methods in unbounded domains using
  {L}aguerre functions.
\newblock SIAM Journal on Numerical Analysis \textbf{38}, 1113--1133 (2001)

\bibitem{shen:2009b}
Shen, J., Tang, T., Wang, L.L.: {S}pectral {M}ethods: {A}lgorithms, {A}nalysis
  and {A}pplications, \emph{Springer {S}eries in {C}omputational
  {M}athematics}, vol.~40.
\newblock Springer (2011)

\bibitem{shen:2009}
Shen, J., Wang, L.L.: Some recent advances on spectral methods for unbounded
  domains.
\newblock Communications in Computational Physics \textbf{5}, 195--241 (2009)

\bibitem{tatari:2014}
Tatari, M., Haghighi, M.: A generalized laguerre--legendre spectral collocation
  method for solving initial-boundary value problems.
\newblock Applied Mathematical Modelling \textbf{38}(4), 1351--1364 (2014)

\bibitem{vismara:2020}
Vismara, F.: A coupled scheme for the solution of parabolic problems on
  unbounded domains.
\newblock Master's thesis, Politecnico di Milano (2020).
\newblock Available at:
  \texttt{https://www.politesi.polimi.it/handle/10589/166694} (last accessed 7
  December 2020)

\bibitem{wang:2009}
Wang, Z.Q., Guo, B.Y., Wu, Y.N.: Pseudospectral method using generalized
  {L}aguerre functions for singular problems on unbounded domains.
\newblock Discrete and Continuous Dynamical Systems Series {B} \textbf{11}(4),
  1019--1038 (2009)

\bibitem{wheeler:1978}
Wheeler, M.: An elliptic collocation-finite element method with interior
  penalties.
\newblock SIAM Journal of Numerical Analysis \textbf{15}, 152--161 (1978)

\bibitem{zhuang:2010}
Zhuang, Q., Shen, J., Xu, C.: A coupled {L}egendre--{L}aguerre
  spectral--element method for the {N}avier--{S}tokes equations in unbounded
  domains.
\newblock Journal of Scientific Computing \textbf{42}(1), 1--22 (2010)

\bibitem{zhuang:2010b}
Zhuang, Q., Xu, C.: Legendre--laguerre coupled spectral element methods for
  second-and fourth-order equations on the half line.
\newblock Journal of computational and applied mathematics \textbf{235}(3),
  615--630 (2010)

\end{thebibliography}

\end{document}